\documentclass[11pt]{amsart}
\newsymbol\subsetneqq 2324
\newtheorem{thm}[equation]{Theorem}
\newtheorem{pro}[equation]{Proposition}
\newtheorem{cor}[equation]{Corollary}
\newtheorem{lem}[equation]{Lemma}
\newtheorem{exa}[equation]{Example}

\newtheorem{DEF}[equation]{Definition}
\newtheorem{rem}[equation]{Remark}

\usepackage{latexsym}

\def\mod{\hbox{mod}}
\def\andd{\quad\hbox{and}\quad}

\def\v{{\mathcal V}}

\def\vd{\dot{\mathcal V}}

\def\fm{(\cdot,\cdot)}
\def\a{\alpha}

\def\w{{\mathcal W}}

\def\sub{\subseteq}
\def\rd{\dot{R}}

\def\rt{\tilde{R}}

\def\adot{\dot{\alpha}}

\def\lam{\lambda}

\def\1k{\frac{1}{k}}

\def\la{\langle}
\def\ra{\rangle}
\def\rds{\dot{R}_{sh}}
\def\rdl{\dot{R}_{lg}}
\def\rs{R_{sh}}

\def\d{\delta}

\def\b{\beta}

\def\qed{\hfill$\Box$}

\def\sg{\sigma}
\def\pia{\Pi(A_{1})}

\def\pib{\Pi(B_{2})}

\def\rtimes{R^{\times}}
\def\gg{{\mathcal G}}
\def\ggi{{\gg_{\bar{i}}}}
\def\ggj{{\gg_{\bar{j}}}}
\def\hh{{\mathcal H}}
\def\hhi{{\hh_{\bar{i}}}}
\def\hhj{{\hh_{\bar{j}}}}
\def\gs{\gg^{\sg}}
\def\hs{\hh^{\sg}}
\def\hhstar{\hh^{\star}}
\def\hsstar{(\hs)^{\star}}
\def\ta{t_\a}
\def\quadd{\quad\quad}
\def\at{\tilde{\a}}
\def\rs{R^{\sg}}
\def\rstimes{(\rs)^{\times}}
\def\rzero{R^{0}}
\def\rszero{(\rs)^{0}}
\def\ad{\hbox{ad}}
\def\pia{\pi(\a)}
\def\pib{\pi(\b)}
\def\pig{\pi(\gamma)}
\def\at{\tilde{\a}}
\def\bt{\tilde{\b}}
\def\pir{\pi(R)}
\def\pirtimes{\pir^{\times}}
\def\dt{\tilde{\delta}}
\def\vpi{\v^{\pi}}
\def\gt{\tilde{\gg}}
\def\vs{\v^{\sg}}
\def\gc{\gg_{c}}
\def\gsc{\gc^{\sg}}
\def\gipia{\gg_{\bar{i},\pia}}

\def\gspia{\gs_{\pia}}
\def\gspib{\gs_{\pib}}
\def\gspig{\gs_{\pig}}
\def\riso{R_{\hbox{iso}}}
\def\rsiso{\rs_{\hbox{iso}}}

\def\rsi{R^{\sg}_{i}}
\def\vsi{\v^{\sg}_{i}}
\def\rsib{\bar{R}^{\sg}_{i}}
\def\vsib{\bar{\v}^{\sg}_{i}}
\def\piib{\bar{\Pi}_{i}}
\def\piid{\dot{\Pi}_{i}}
\def\vsid{\dot{\v}^{\sg}_{i}}
\def\rsid{\dot{R}^{\sg}_{i}}
\def\rniso{R_{\hbox{niso}}}
\def\rsniso{\rs_{\hbox{niso}}}

\def\rsubt{R_{\hbox{t}}}

\def\dj{{\mathcal{D}}_{j}}
\def\di{{\mathcal D}_{i}}
\def\ri{R_{i}}
\def\hsi{\hs_{i}}

\def\ti{{\mathcal T}_{i}}
\def\ki{{\mathcal K}_{i}}
\def\gsat{\gs_{\at}}
\def\mi{{\mathcal M}_{i}}
\def\gsi{\gs_{i}}

\def\i{{\mathcal I}}
\def\be{{\bf e}}
\def\bq{{\bf q}}
\def\bbbc{{\mathbb C}}
\def\bbbz{{\mathbb Z}}
\def\bbbr{{\mathbb R}}
\def\zn{{\mathbb Z}^{\nu}}
\def\aa{\mathcal A}

\def\zeq{Z_{\be,\bq}}
\def\kk{\mathcal K}
\def\mna{M_{n}(\aa)}
\def\hd{\dot{\hh}}
\def\slna{sl_n(\aa)}
\def\tr{\hbox{tr}}
\def\cc{{\mathcal C}}
\def\dd{\mathcal D}
\def\ep{\epsilon}
\def\ks{\kk^{\sg}}
\def\mna{M_{n}(\aa)}
\def\hds{\hd^{\sg}}
\def\xd{x^{\d}}

\def\kd{\kk^{\d}}
\def\aff{\hbox{Aff}}
\def\hht{\tilde{\mathcal H}}
\def\ati{\at_i}
\def\gamt{\tilde{\gamma}}
\def\vsd{\dot{\v}^{\sg}}

\begin{document}
\markboth{FIXED POINT SUBALGEBRAS OF EALA'S}{S. AZAM, S. BERMAN,
M. YOUSOFZADEH}

\centerline{\bf FIXED POINT SUBALGEBRAS OF EXTENDED}
\centerline{\bf AFFINE LIE ALGEBRAS}

\vspace{.5cm}\centerline{SAEID AZAM\footnote[1]{Research number
810821, University of Isfahan, Iran.}$^{,}$\footnote[2]{This
research was in part supported by a grant from IPM.}}
\centerline{STEPHEN BERMAN\footnote[3]{The author gratefully
acknowledge the support of the National Sciences and Engineering
Research Council of Canada}}\centerline{MALIHE YOUSOFZADEH}

%\footnote{1991 Mathematics Subject classification. Primary 17B65,
%Secondary 17B67, 17B40.}
%\footnotetext[1]{Research number
%810821, University of Isfahan, Iran.} \footnote[2]{This research
%was in part supported by a grant from IPM.}
%\date{}
%\maketitle
%\begin{abstract}

\vspace{1cm} \noindent ABSTRACT. It is a well known result that
the fixed point subalgebra of a finite dimensional complex simple
Lie algebra under a finite order automorphism is a reductive Lie
algebra so is a direct sum of finite dimensional simple Lie
subalgebras and an abelian subalgebra. We consider this for the
class of extended affine Lie algebras and are able to show  that
the fixed point subalgebra of an extended affine Lie algebra
under a finite order automorphism (which satisfies certain
natural properties) is a sum of extended affine Lie algebras (up
to existence of some isolated root spaces), an abelian subalgebra
and a subspace which is contained in the centralizer of the core.

\setcounter{section}{-1}
\section{Introduction}
In 1955, A. Borel and G.D. Mostow [BM] proved that the fixed
point subalgebra of a finite dimensional complex simple Lie
algebra under a finite order automorphism is a reductive Lie
algebra. A natural question which arises here is  what  can we
say about the fixed points of a finite order automorphism of an
{\it extended affine Lie algebra} (EALA for short). EALA's are
natural generalizations of finite dimensional complex simple Lie
algebras and affine Kac--Moody Lie algebras. They are
axiomatically defined (see [AABGP], [HK-T]) and the axioms
guarantee the existence of analogues of Cartan subalgebras, root
systems, invariant forms, etc. A root of an EALA is called {\it
isotropic} if it is orthogonal to itself, with respect to the
form. The dimension of the real span of the isotropic roots is
called the {\it nullity} of the Lie algebra. A finite dimensional
simple Lie algebra is an EALA of nullity zero, and an EALA is an
affine Lie algebra if and only if its nullity is $1$ (see [ABGP]
for details). Thus EALAs form a natural class of algebras in
which to consider extensions of the result of Borel and Mostow.

Here we would like to explain a procedure which has been the most
general theme of constructing affine Lie algebras and their
generalizations, since the birth of Kac--Moody Lie algebras in
1968.

Let $(\gg,\fm,\hh)$ be an EALA with root system $R$ (in
particular $\gg$ can be a finite dimensional simple Lie algebra
or an affine Lie algebra). Let $\sg$ be a finite order
automorphism of $\gg$ which stabilizes $\hh$ and leaves the form
invariant. Assume also that the fixed point subalgebra of $\hh$
(with respect to $\sg$) is a Cartan subalgebra of the fixed point
subalgebra of $\gg$. Consider the Lie algebra
$$
\aff(\gg):=(\gg\otimes \bbbc [t,t^{-1}])\oplus\bbbc c\oplus\bbbc
d,
$$
where $c$ is central, $d=t\frac{d}{dt}$ is the degree derivation
so that $[d,x\otimes t^n]=nx\otimes t^n$, and  multiplication is
given by $$ [x\otimes t^n,y\otimes t^m]=[x,y]\otimes
t^{n+m}+n(x,y)\d_{m+n,0}c. $$ Extend the form $\fm$ to
$\aff(\gg)$ so that $c$ and $d$ are naturally paired. Then the
triple
$$
(\aff(\gg),\fm,\hh\oplus\bbbc c\oplus\bbbc d)
$$
is again an EALA with root system $\rt=R+\bbbz\d$ where $\d$ is
the linear functional on $\hh\oplus\bbbc c\oplus\bbbc d$ defined
by $\d(d)=1$ and $\d(\hh\oplus\bbbc c)=0$. Extend $\sg$ to an
automorphism of $\aff(\gg)$ by
$$
\sg(x\otimes t^i+ rc+sd)=\zeta^{-i}\sg(x)\otimes t^i+rc+sd,
$$
where $\zeta=e^{2\pi\sqrt{-1}/m}$ and $\sg^{m}=\hbox{id}$. Then
the fixed point subalgebra of $\aff(\gg)$, under $\sg$, is an
EALA (see [ABP] and Example \ref{3-12}). As one can see in [K] and
[H], when $\gg$ is a finite dimensional simple Lie algebra, all
affine Kac--Moody Lie algebras can be constructed this way. In
fact all the examples of EALA's which are presented in [W], [P]
and [HK-T] can be obtained from the above procedure. It is also
shown in [A2] that most of the examples of EALA's constructed in
[AABGP] can be put in the above context (see Section 3 and [A2]).

In this paper, we consider the theorem of [BM] for the class of
EALA's. Namely we investigate what is the fixed point subalgebra
of an EALA under a finite order automorphism which satisfies
certain conditions. In a very real sense this paper found it's
inspiration from the paper [ABP] which also studies Lie algebras
constructed from finite order automorphisms of an EALA. The
difference is this; in [ABP] the algebras studied are
affinizations of the EALA one starts with while here we study
fixed point subalgebras. But our basic techniques are similar to
[ABP] and, as the reader will see, we will rely on some of the
results from [ABP] in certain crucial places. Thus, although both
studies are independent of each other there is certainly an
interdependence.

 In Section 1 of this paper, we state a modified version of
the definition of an extended affine root system (EARS for short).
The usual definition has indecomposability built into it. We need
to consider root systems which are unions of such root systems so
have broken down the usual definition. This should not cause the
reader any difficulties. The definition given here is basically
the same as [AABGP] except for a rearrangement of the axioms.
Part of the reason this is done is to make more compatible the
two different definitions of EARSs  found in the literature. It
is shown that an EARS, modulo some isotropic roots, is a union of
a finite number of irreducible EARS's which are orthogonal with
respect to the form.

In Section 2, which forms the core of the paper, we show that the
fixed point subalgebra $\gs$ of an EALA $\gg$ under a finite
order automorphism $\sg$ which satisfies certain conditions is of
the form $\gs=\sum_{i=1}^{k}\gsi\oplus\w\oplus\i$, where for each
$i$, $\gsi$ satisfies axioms EA1--EA5(a) of an EALA, $\i$ is a
subspace of $\gs$ contained in the centralizer of the core of
$\gs$, and $\w$ is an abelian subalgebra contained in the
centralizer of $\gsi$ for each $i$ (see Theorem \ref{2-27}). This
agrees with the result of [BM] when $\gg$ is a finite dimensional
simple Lie algebra over $\bbbc$, as in this case $\i=\{0\}$ and
for each $i$, $\gsi$ is a finite dimensional simple Lie algebra
(see Corollary \ref{2-29}). In Lemma \ref{2-16b} certain relations
between the core of $\gg$ and the core of $\gs$ and some results
regarding the tameness of $\gs$ in terms of the tameness of $\gg$
are obtained.

In Section 3, the last section, a large number of examples are
presented. Examples \ref{3-2}-\ref{3-5} illustrate how the terms
$\gsi$'s, $\w$ and $\i$ (see Theorem \ref{2-27}) appear as the
fixed points of automorphisms. In \ref{3-6}-\ref{3-9} and
\ref{3-11} some examples from [A2] are restated in such a way
that they fit in the setting presented in Section 2. In
\ref{3-12}, an example regarding the results in [ABP] is provided.
This example  shows how the affinizations in [ABP] can be viewed
as fixed point subalgebras of the type considered here.

The authors would like to thank Professor B. Allison and Professor
A. Pianzola for some helpful discussions. It was at the Fields
Institute Program on Infinite Dimensional Lie Theory and Its
Applications where this work was originally conceived
\medskip

\section{Terminology and prerequisites}
There are two different definitions for the term {\it extended
affine root systems} (EARS for short) in the literature. This
term first was introduced by K. Saito [Sa1] in 1985. In 1997,
another definition for EARS was introduced in [AABGP]. In [A1],
to distinguish between these two definitions,the author used
SEARS to refer to the definition
 introduced by K. Saito and used EARS for the one
which is introduced in [AABGP]. The basic difference between
SEARS and EARS is that in SEARS all roots are nonisotropic while
in EARS there are isotropic roots as well. Also in [AABGP] it is
assumed that an EARS is both reduced and indecomposable while
this is not assumed in the case of a SEARS. In [A1], it is shown
that there is a one to one correspondence between reduced
indecomposable SEARS and EARS. In fact, if we drop the axioms
which are related to an EARS being reduced and indecomposable
(axioms R4 and R7 from [AABGP])  it follows from [A1] that there
is a one to one correspondence between SEARS and EARS. As we will
see in the sequel, such root systems will naturally arise as the
root systems of the fixed point subalgebras (under certain finite
order automorphisms) of some extended affine Lie algebras. Also,
such root systems arise in [Y].
%A study of EARS (in the sense of Definition \ref{EARS}) is under
%progress by [AKY].

We start by modifying the definition of an EARS. Our new
definition is basically the same as [AABGP], except, since we
want to consider more general types of root systems, we leave off
the axioms concerning reducibility and indecomposability.

\begin{DEF}\label{EARS}
{\rm Let $\v$ be a finite dimensional real vector space with a
nontrivial positive semidefinite symmetric bilinear form $(.,.)$
and let $R$ be a subset of $\v$. Let
$$
R^{\times}=\{\a\in R:(\a,\a)\not=0\}\quad\hbox{and} \quad
R^{0}=\{\a\in R:(\a,\a)=0\}.
$$
Then $R=R^{\times}\uplus R^{0}$ where $\uplus$ means disjoint
union. Then we will say $R$ is an {\it extended affine root
system} (EARS) in $\v$ if $R$ satisfies the following 4 axioms:

\noindent (R1) $R=-R$,

\noindent (R2)$\;R \hbox{ spans }\v$,

\noindent (R3)$\;$ $R$ is discrete in $\v$,

\noindent (R4)$\;$ if  $\a\in R^{\times}$ and $\b\in R$, then
there exist $d, u\in {\mathbb Z}_{\geq 0}$ such that
$$\{ \b+n\a:n \in {\bbbz} \} \cap R = \{ \b-d\a,\ldots,\b+u\a \}
\hbox{ and }d-u=2\frac{(\a,\b)}{(\a,\a)}.
$$

The EARS $R$ is called {\it tame}  if it satisfies:

\noindent (R5)$\;$ for any $\d\in R^{0}$, there exists $\a\in
R^{\times}$ such that $\a+\d\in R$. We say a root satisfying this
condition is
 {\it nonisolated} and call isotropic roots which do not satisfy this
{\it isolated}.

The EARS $R$ is called {\it indecomposable} if it satisfies:

\noindent (R6)$\;$ $R^{\times}$ cannot be decomposed into a
disjoint union of two nonempty subsets which are orthogonal with
respect to the form.

A tame indecomposable EARS $R$ is called {\it irreducible}.
Finally, the EARS $R$ is called {\it reduced} it satisfies:

\noindent (R7)$\;\a\in R^{\times} \Rightarrow 2\a\not\in R$.}
\end{DEF}

Since the form is nontrivial, it follows from R2 that
$\rtimes\not=\emptyset$. This, together with R1 and R4, implies
that $0\in R$. Note that in an EARS as used here could have  both
isolated and nonisolated roots. Thus Axiom R5 indicates that in a
tame EARS isotropic roots are nonisolated. By [A1], there is a
one to one correspondence between irreducible (reduced) EARS and
indecomposable (reduced) SEARS.

\begin{lem}\label{yos}
Let $R$ be an EARS and $R_1$ be a subset of $R$ with
$R_{1}^{\times}\not=\emptyset$. Suppose that

(a) $R_1=-R_1$,

(b) $\{\d\in R^0\mid\a'+\d\in R_1\hbox{ for some }\a'\in
R_{1}^{\times}\}\sub R_1$,

(c) $\a'\in R_1$, $\b\in R$, $(\a',\b)\not=0\Longrightarrow \b\in
R_1$.

Then $R_1$ is an EARS in its real span. Moreover, if we set
$$
R_{1}'=R_{1}^{\times}\cup(\la R_1\ra\cap R^0),
$$
then $R_{1}'$ is also an EARS in the real span of $R_1$ ($\la
R_1\ra$ denotes the ${\bbbz}$-span of $R_1$).
\end{lem}

\noindent Proof. Clearly (R1)--(R3) hold for $R_1$. We now check
$R4$. Let $\a'\in R_{1}^{\times}$ and $\b'\in R_1$. Since $R4$
holds for $R$, it is enough to show that for $n\in{\bbbz}$,
$$
\b'+n\a'\in R\Longrightarrow\b'+n\a'\in R_1.
$$
Since $\b'\in R_1$, we may assume that $n\not=0$, and by (a), we
may also assume that $n>0$. So let $\b'+n\a'\in R$, $n>0$. If
$\b'+n\a'\in\rtimes$, then $(\b'+n\a',\b')\not=0$ or
$(\b'+n\a',\a')\not=0$. In either cases, we get from (c) that
$\b'+n\a'\in R_1$. Next, let $\b'+n\a'\in R^0$. Since $R4$ holds
for $R$ and $n>0$, we have $\b'+(n-1)\a'\in\rtimes$. So repeating
our previous argument we get $\b'+(n-1)\a'\in R_1$. Since
$$
\b'+n\a'+(-\a')=\b'+(n-1)\a'\in R_{1}^{\times}
$$
it follows from (a) and (b) that $\b'+n\a'\in R_1$. This
completes the proof of the first assertion.

Next let $R_{1}'$ be as in the statement. Clearly $R_1$ and
$R_{1}'$ have the same real span. Since
$R_{1}^{\times}=(R_{1}')^{\times}$, it is easy to check that
$R_{1}'$ satisfies conditions (a)-(c), and so is an EARS.
\hfill\qed \vspace{3mm}

Let $R$ be an EARS in $\v$. Let $\v^0$ be the radical of the form
$\fm$ on $\v$. Set $\bar{\v}:=\v/{\v^{0}}$, and let $\bar{v}$ be
the image of an element $v$ of $\v$ under the projection map
$\v\longrightarrow\bar{\v}$. For $\a$, $\b\in\v$ define
$(\bar{\a},\bar{\b}):=(\a,\b)$. Then $\fm$ is positive definite
on $\bar{\v}$, and by [AABGP, Chapter II] $\bar{R}$ is a finite
root system in $\bar{\v}$. Note that it is not assumed here that
$R$ is irreducible , but the same argument as in [AABGP, Chapter
II] shows that $\bar{R}$ is a finite root system (which is not
necessarily irreducible).

For an EARS $R$, we set
\begin{equation}\label{riso}
\begin{array}{l} \riso=\{\d\in
R^0\mid\a+\d\not\in R\hbox{ for any }
\a\in\rtimes\},\\
 \rniso=\{\d\in R^0\mid\a+\d\in R\hbox{ for some }
 \a\in\rtimes\}=R^0\setminus\riso.\\
\end{array}
\end{equation}
That is $\riso$ ($\rniso$) is the set of isolated (nonisolated)
isotropic roots of $R$, so $R$ is tame if and only if
$\riso=\emptyset$. We also set
$$\rsubt=\rtimes\cup\rniso.$$
Then $R_{t}^{\times}=\rtimes$ and
$$
R=\rtimes\uplus\rniso\uplus\riso=\rsubt\uplus\riso.
$$

\begin{lem}\label{indecom}
Let $R$ be an EARS. Then

(i) $\rsubt$ is a tame EARS in its real span.

 (ii) $R=(\cup_{i=1}^{k}R_{i})\cup\riso$,
where each $R_i$ is an irreducible EARS, and for $i\not=j$,
$R_{i}$ and $R_{j}$ are orthogonal with respect to the form.
Furthermore, if we set
\begin{equation}\label{EARSa}
R'_{i}=\ri^{\times}\uplus(\la\ri\ra\cap R^{0}),
\end{equation}
then $R'_{i}$ is an indecomposable EARS.
\end{lem}

\noindent {\bf Proof.} (i) It follows from definition of $R_t$
that conditions (a)-(c) of Lemma \ref{yos} hold for $R_{t}$. So
$R_t$ is an EARS. It is clear from definition that $R_t$ is tame.
%It is clear that axioms R1-R3 hold for $\rsubt$. Since R4 holds
%for $R$, to show that R4 holds for $R_{t}$, it is enough to show
%that if $\a\in R_{t}^{\times}=\rtimes$ and $\b\in\rsubt$ then
%$$
%\{\b+n\a\mid n\in\bbbz\}\cap R=\{\b+n\a\mid
%n\in\bbbz\}\cap\rsubt.
%$$
%This is clear that the right hand side is a subset of the left
%hand side. To see the reverse inclusion, let $\b+n\a\in R$. If
%$\b+n\a$ is nonisotropic then $\b+n\a\in\rtimes\sub\rsubt$. If
%$n=0$, then $\b+n\a=\b\in \rsubt$. So we may assume that $n>0$ and
%$\b+n\a\in\rzero$. Since R4 holds for $R$, we have $\b+(n-1)\a\in
%R$. So the isotropic root $\b+n\a$ is attached to the nonisotropic
%root $-\a$ and is nonisolated. That is
%$\b+n\a\in\rniso\sub\rsubt$. From definition it is clear that
%$\rsubt$ is tame.%

(ii) We have $R=\rsubt\uplus\riso$, and by part (i), $\rsubt$ is
a tame EARS. Since $\bar{\rsubt}=\bar{R}$ is a finite root system,
we have $\bar{\rsubt}\setminus\{0\}=\uplus_{i=1}^{k}\bar{R}_{i}$,
where for each $i$, $\bar{R}_{i}\cup\{0\}$ is an irreducible
finite root system, and $\bar{R}_{i}$'s are orthogonal with
respect to the form. Let $R^{\times}_i$ be the preimage of
$\bar{R}_{i}$ under the projection map ${\bar{\;}}$. Then
\begin{equation}\label{1}
\rsubt^{\times}=\uplus_{i=1}^{k} \rtimes_{i}\quad\hbox{and}\quad
(\rtimes_i,\rtimes_{j})=\{0\} \hbox{ if } i\not= j.
\end{equation}
Set
\begin{equation}\label{2}
R_i=\rtimes_i\cup\{\d\in R_{t}^{0}| \d+\a\in R_{t}\hbox{ for some
} \a\in R^{\times}_i\}.
\end{equation}
Since $R_{t}$ is tame, the isotropic roots of $R_{t}$ are
nonisolated so
$$R_{t}=\cup_{i=1}^{k}R_i.
$$
It is easy to see that conditions (a)-(c) of Lemma \ref{yos} hold
for $R_i$. Thus $R_i$ is an EARS.
% Let $\v_i$ be the real span of
%$R_i$. We want to show that $R_i$ is an irreducible EARS in
%$\v_i$. It is clear that axioms R1-R3 and R5-R6 hold for $R_i$.
%So it only remains to check R4. Let $\a\in\rtimes_{i}$ and $\b\in
%R_i$. Then R4 holds for $\a$ and $\b$ as elements of $R$.
%Therefore it is enough to show that
%\begin{equation}\label{string}
%\b+k\a\in R\Longrightarrow\b+k\a\in R_i\quad\quad(k\in{\mathbb Z}).
%\end{equation}
%Let $\b+k\a\in R$. As $\b\in R_i$, we may assume that $k\not=0$.
%If $(\a,\b)=0$, then
%$$
%(\b+k\a,\a)=k(\a,\a)\not=0.
%$$
%So by (\ref{1}), $\b+k\a\in R_i$. Next assume $(\a,\b)\not=0$. By
%replacing $\a$ with $-\a$ if necessary, we may assume that $k>0$.
%It is enough to show that $\b+\a\in R_i$, since if $k>1$, we can
%use an inductive argument to conclude that $\a+j\b\in R_i$ for
%$1\leq j\leq k$. Now if
%$$
%(\b+\a,\a)\not=0,\quad\hbox{or}\quad(\b+\a,\b)\not=0,
%$$
%then by (\ref{1}), $\b+\a\in R_i$, and we are done. Otherwise,
%$(\b+\a,\b+\a)=0$, so $\b\in \rtimes_i$ and $\b+\a\in R^{0}$. Now
%$-\a\in \rtimes_i$ and
%$$
%-\a+(\a+\b)=\b\in\rtimes_i.
%$$
%So by the way $R_i$ is defined (see (\ref{2})), $\a+\b\in R_i$.%
By Lemma \ref{yos} the last assertion also holds.
%For the last
%assertion, first note that  $R'_{i}$ and $\ri$ have the same
%${\mathbb Z}$--span, so R1--R3 and R6 hold for $R'_{i}$ as they
%hold for $\ri$. We only need to show that R4 holds. So let $\a\in
%(R'_{i})^{\times}=\ri^{\times}$, and $\b\in R'_{i}$. As elements
%of $R$, R4 holds for $\a$ and $\b$. Consider an element $\b+n\a$
%in the $\a$--string through $\b$. As $\b\in R'_{i}$, we may
%assume that $n\not=0$. If $\b+n\a$ is isotropic, then it is in
%the intersection of the ${\mathbb Z}$--span of $\ri$ with
%$R^{0}$, so $\b+n\a\in R'_{i}$, by the way $R'_{i}$ is defined.
%Next suppose that $\b+n\a$ is nonisotropic. If $\b$ is
%nonisotropic, then $\a,\b\in\ri$ and so $\b+n\a\in \ri$ as R4
%holds for $\ri$. But $\ri\sub R'_{i}$. Finally, let $\b$ be
%isotropic. Then $(\b+n\a,\a)\not=0$ so $\b+n\a\in\ri\sub R'_{i}$.
%This finishes the proof.
\hfill\qed\vspace{5mm}

\section{Fixed Point Subalgebras}
In this section we study the structure of fixed point subalgebra
of an extended affine Lie algebra under a finite order
automorphism which satisfies certain conditions. For a systematic
study of extended affine Lie algebras reader is referred to
[AABGP].

An {\it extended affine Lie algebra} (EALA for short) is a triple
$(\gg,\fm,\hh)$ where $\gg$ and $\hh$ are two Lie algebras over
$\bbbc$ and $\fm$ is a complex valued bilinear form on $\gg$
satisfying the following five axioms:

EA1. The form $\fm$ is symmetric, nondegenerate and invariant on
$\gg$.

EA2. $\hh$ is a nontrivial finite dimensional abelian subalgebra
of $\gg$ which is self-centralizing and $\ad(h)$ is diagonalizable
for all $h\in\hh$.

Consider the root space decomposition $\gg=\oplus_{\a\in
\hh^\star}\gg_{\a}$. Then
$R=\{\a\in\hh^\star\mid\gg_\a\not=\{0\}\}$ is called the root
system of $\gg$. Let $\rtimes$ and $R^0$ be as in Definition
\ref{EARS}. The next three axioms are as follows:

EA3. For any $\a\in\hh^\star$ and $x\in\gg_{\a}$, $\ad_{\gg}(x)$
is locally nilpotent on $\gg$.

EA4. $R$ is a discrete subset of $\hh^\star$.

EA5a. $\rtimes$ is indecomposable (in the sense of Definition
\ref{EARS}),

EA5b. Isotropic roots of $R$ are nonisolated (in the sense of
Definition \ref{EARS}).

Consider a fixed EALA $(\gg,\fm,\hh)$ with root system $R$ and the
corresponding root space decomposition
$\gg=\oplus_{\a\in\hh^\star}\gg_\a$. It is shown in [AABGP,
Chapter I] that $R$ is a reduced irreducible EARS (in the sense
of Definition \ref{EARS}). Consider a fix automorphism $\sg$ of
$\gg$ and set
\begin{equation}\label{2-1}
%\begin{array}{l}
\gs=\{x\in\gg|\sg(x)=x\}\andd \hs=\{h\in\hh|\sg(h)=h\}.
%\end{array}
\end{equation}
That is $\gs$ (resp. $\hs$) is the fixed point subalgebra of
$\gg$ (resp. $\hh$) with respect to $\sg$.

Let $m\geq 1$ and suppose that

A1.$\quad\sg^m=1$.

A2.$\quad\sg(\hh)=\hh$.

A3.$\quad(\sg(x),\sg(y))=(x,y)\quad$ for all $x,y\in\gg$.

A4.$\quad C_{\gs}(\hs)=\hs$.

We start by recording some facts related to A1-A3, so only need
to assume A1-A3 hold for now. Let $\bar{i}$ denote the image of
$i\in {\mathbb Z}$ in ${\mathbb Z}/m{\mathbb Z}$. From A1 and A2,
We have
\begin{equation}\label{2-2}
\gg=\oplus_{i=0}^{m-1}\ggi\andd\hh=\oplus_{i=0}^{m-1}\hhi,
\end{equation}
where $\ggi$ (resp $\hhi$) is the eigenspace corresponding to the
i-th power of the m-th root of unity $\zeta=e^{2\sqrt{-1}\pi/m}$.
Then
$$
\ggi=\{x\in\gg|\sg(x)=\zeta^{i}x\}.
$$
Note that $\gs=\gg_{\bar{0}}$ and $\hs=\hh_{\bar{0}}$. It follows
from A3 that
\begin{equation}\label{2-3}
(\ggi,\ggj)=\{0\}=(\hhi,\hhj)\quad\hbox{ if
}\bar{i}+\bar{j}\not=0.
\end{equation}
Set $\gg^{c}=\oplus_{i=1}^{m-1}\ggi$ and
$\hh^{c}=\oplus_{i=1}^{m-1}\hhi$. Then $\gg=\gs\oplus\gg^{c}$ and
$\hh=\hs\oplus\hh^{c}$. Moreover,
\begin{equation}\label{2-4}
(\gs,\gg^{c})=\{0\}=(\hs,\hh^{c}).
\end{equation}

Note also that $\sg$ induces an automorphism
$\sg\in\hbox{Aut}(\hhstar)$ by
\begin{equation}\label{2-4a}
\sg(\a)(h)=\a(\sg^{-1}(h)),\quad\hbox{for }\a\in\hhstar\hbox{ and
}h\in\hh.
\end{equation}
Let $\pi$ denote the projection map from $\gg$ (resp. $\hh$) onto
$\gs$ (resp. $\hs$). Then for $x\in\gg$ we have
$x-\pi(x)\in\gg^{c}$. Since $\gs$ is stable under $\sg$ and
$\sg\pi=\pi$, we obtain
$$
\sum_{i=0}^{m-1}\sg^{i}(x-\pi(x))\in\gg^{c}\cap\gs=\{0\}.
$$
So
\begin{equation}\label{2-4b}
\pi(x)=\frac{1}{m}\sum_{i=0}^{m-1}\sg^{i}(x),\quadd(x\in\gg).
\end{equation}
Similarly,
\begin{equation}\label{2-5}
\pi(h)=\frac{1}{m}\sum_{i=0}^{m-1}\sg^{i}(h),\quad\quad(h\in\hh).
\end{equation}
The map $\pi$ induces a map, denoted again by $\pi$, from the dual
space $\hhstar$ of $\hh$ onto the dual space $\hsstar$ of $\hs$.
Namely, for $\a\in\hhstar$ and $h\in\hs$, define
\begin{equation}\label{2-6}
\pi(\a)(h)=\a(\pi(h)),
\end{equation}
that is $\pi(\a)$ for $\a\in\hh^\star$ is the restriction of $\a$
to $\hs$. Then from (\ref{2-4a}) we have, for $\a\in\hhstar$ and
$h\in\hs$,
$$
\pi(\a)(h)=\a(\pi(h))=\a(\frac{1}{m}\sum_{i=0}^{m-1}\sg^{i}(h))
=\frac{1}{m}\sum_{i=0}^{m-1}\sg^{i}(\a)(h).
$$
Thus for $\a\in\hhstar$, we have
\begin{equation}\label{2-7}
\pi(\a)=\frac{1}{m}\sum_{i=0}^{m-1}\sg^{i}(\a).
\end{equation}

Since both $\gg$ and $\gg_{\bar{i}}$ are $\hs$--modules, we have
\begin{equation}\label{2-7a}
\gg=\sum_{\at\in(\hs)^{\star}}\gg_{\at}=\sum_{\a\in\hh^{\star}}\gg_{\pia}=\sum_{\a\in
R}\gg_{\pia}\andd\gg_{\bar{i}}=\sum_{\a\in R}\gg_{\bar{i},\pia},
\end{equation}
where
$$ \gg_{\at}=\{x\in\gg\mid [h,x]=\at(h)x\quad\hbox{for all
}h\in\hs\}$$ and
$$\gg_{\bar{i},\pia}=\gg_{\bar{i}}\cap\gg_{\pia}.
$$

Note that $\gs$ as an $\hs$-submodule of $\gg$ has a weight space
decomposition
\begin{equation}\label{2-12}
\gs=\bigoplus_{\at\in\hsstar}\gs_{\at},
\end{equation}
where
$$\gs_{\at}=\gs\cap\gg_{\at}.
$$

Set
$$
\rs=\{\at\in\hsstar|\gs_{\at}\not=\{0\}\}.
$$
Then
$$\rs\sub\pi(R).
$$
As we will see in the next section, in many examples $\rs$ is in
fact a proper subset of $\pi(R)$. Denote the set of nonisotropic
(isotropic) roots of $\rs$ by $\rstimes$ ($\rszero$),
respectively. That is,
$$
\rstimes=\{\at\in\rs |(\at,\at)\not=0\}
\andd\rszero=\{\at\in\rs|(\at,\at)=0\}.
$$

It is easy to see that
\begin{equation}\label{2-10}
\sg(\gg_{\a})=\gg_{\sg(\a)}\quadd(\a\in R),
\end{equation}
and so
$$
\sg(R)=R.
$$
Also we have
$$
\pi(\gg_{\a})\sub\gs_{\pi(\a)}.
$$
Thus
\begin{equation}\label{2-11}
\gs_{\pi(\a)}=\sum_{\stackrel{\b\in
R}{\pi(\b)=\pi(\a)}}\pi(\gg_{\b}).
\end{equation}

Since the form $\fm$ is nondegenerate (by (\ref{2-4})) and
invariant on $\gs$, it follows that
\begin{equation}\label{2-13z}
(\gs_{\at},\gs_{\bt})=\{0\}\quad\quad\hbox{unless}\quad\at+\bt=0.
\end{equation}
In particular, we have
\begin{equation}\label{2-13zz}
\rs=-\rs.
\end{equation}
Also for any $\a,\b\in R$,
\begin{equation}\label{2-13a}
[\gspia,\gspib]=0\quad\hbox{or}\quad[\gspia,\gspib]\sub\gspig,
\end{equation}
for some $\gamma\in R$ with $\pia+\pib=\pig$. Then from
(\ref{2-11}) and (\ref{2-12}), we have
\begin{equation}\label{2-13}
\gs=\sum_{\a\in
R}\gs_{\pi(\a)}=\bigoplus_{\at\in\hsstar}\gs_{\at}=
\bigoplus_{\at\in\rs}\gs_{\at}.
\end{equation}

For $\a\in\hhstar$ let $\ta$ be the unique element in $\hh$ which
represents $\a$ via the nondegenerate form $\fm$. That is,
$$
\a(h)=(\ta,h)\quadd(h\in\hh).
$$
Then for $\a\in\hh^{\star}$,
\begin{equation}\label{2-8}
\sg(\ta)=t_{\sg(\a)}\andd \pi(t_\a)=t_{\pia}.
\end{equation}
So
\begin{equation}\label{2-9}
\hs=\{\ta|\a\in\hsstar\}\andd\hh^{c}=\{\ta|\a\in(\hh^{c})^{\star}\}.
\end{equation}
Transfer the form $\fm$ to $\hhstar$ by
$$
(\a,\b)=(\ta,t_{\b}),\quadd(\a,\b\in\hhstar).
$$
Now similar to (\ref{2-4}), we have
\begin{equation}\label{2-10-1}
\hhstar=(\hhstar)^{\sg}\oplus(\hhstar)^{c},\andd
((\hhstar)^{\sg},(\hhstar)^{c})=\{0\}.
\end{equation}
Identifying $(\hhstar)^{\sg}$ with $\hsstar$ and $(\hhstar)^{c}$
with $(\hh^{c})^{\star}$, we see that the map $\pi$ on $\hhstar$
is in fact the restriction to $\hsstar$.

\begin{lem}\label{2-14}
Let $\a\in R$ and $(\pi(\a),\pi(\a))\not=0$. Then for
$x\in\gs_{\pi(\a)}$, $\ad_{\gs}(x)$ acts locally nilpotently on
$\gs$.
\end{lem}

\noindent{\bf Proof.} Fix $0\not=x\in\gs_{\pi(\a)}$. By
(\ref{2-13}), it is enough to show that if $\b\in R$ and
$y\in\gs_{\pi(\b)}$, then $(\ad x)^{N}(y)=0$ for sufficiently
large $N$. By (\ref{2-8}) and (\ref{2-9}), $h:=t_{\pia}\in\hs$ and
so
$$
[h,(\ad x)^{k}(y)]=\big(k(\pi(\a),\pi(\a))+\pi(\b)(h)\big)(\ad
x)^{k}(y).
$$
According to (\ref{2-13}) the eigenvalues of $\ad_{\gs}h$ are of
the form $\pi(\gamma)(h)$, $\gamma\in R$. Now
$$
\pi(\gamma)(h)=(\pia,\pi(\gamma))=
\frac{1}{m}\sum_{i=0}^{m-1}(\a,\sg^{i}(\gamma)).
$$
But $\bar{R}$ is a finite root system, so $\{(\a,\gamma)|\gamma\in
R\}=\{(\bar{\a},\bar{\gamma})\mid\gamma\in R\}$ is a finite set.
This shows that $\ad_{\gs}h$ has only a finite number of
eigenvalues, and so $(\ad x)^{N}(y)=0$ for sufficiently large
$N$.\hfill\qed\vspace{3mm}

\begin{lem}\label{2-15}
$\pi(R)$ is a discrete subset of $\hsstar$. In particular $\rs$
as a subset of $\pi(R)$ is discrete.
\end{lem}

\noindent{\bf Proof.} This is clear as
$\pi(\a)=\frac{1}{m}\sum_{i=0}^{m-1}\sg^{i}(\a) \in\frac{1}{m}\la
R\ra$, and  $\la R\ra$ is a discrete subset of
$\hhstar$.\hfill\qed\vspace{3mm}

The parts (i) and (ii) of the following lemma are proved in [ABP].
For part (iii) we need to recall that the structure of the root
system $R$ of type $X_{\ell}$ is of the form
\begin{equation}\label{2-15a}
\begin{array}{ll}
(S+S)\cup(\rd+S)&\hbox{if }X_{\ell}=A_{\ell},\; D_{\ell},\; E_{\ell}\\
(S+S)\cup (\rds+S)\cup(\rdl+L)&\hbox{if }X_{\ell}
=B_{\ell},\;C_{\ell},\;F_4,\;G_2\\
(S+S)\cup (\rds+S)\cup (\rd_{ex}+E)&\hbox{if }X_{\ell}=BC_{1}\\
(S+S)\cup (\rds+S)\cup(\rdl+L)\cup (\rd_{ex}+E)&\hbox{if
}X_{\ell}=BC_{\ell}
\end{array}
\end{equation}
where $\rds$, $\rdl$ and $\rd_{ex}$ are the sets of short, long
and extra long roots of a finite root system $\rd$, and $S$, $L$,
$E$ are the so called semilattices or translated semilattices
involved in the structure of $R$ which interact in a prescribed
way. Moreover, $\rds\cup\rdl\subseteq R$ (for details see [AABGP,
Chapter II].)

\begin{lem}\label{2-16}
Suppose $\pirtimes\not=\emptyset$. Then

(i) For $\b\in \rtimes$ there exists $\a\in R$ such that
$(\pia,\pia)\not=0$ and $(\a,\b)\not=0$.

(ii) $\pirtimes$ is indecomposable. That is, $\pirtimes$ cannot
be written as a disjoint union of two nonempty subsets which are
orthogonal with respect to the form.

(iii) Let $x\in R$ and $(\pi(x),\pi(x))=0$. Then there exists
$\bt\in\pirtimes$ such that $\bt+\pi(x)\in\pirtimes$.
\end{lem}

\noindent{\bf Proof.} For parts (i) and (ii) see [ABP Lemma
3.49]. To see (iii), we consider the two cases (a) $x\in R^{0}$
and $x\in\rtimes$, separately. Let $\dt:=\pi(x)$.

(a) $x\in\rzero=S+S$. Then $x=\d_1+\d_2$ for some $\d_1,\d_2\in
S\sub\rzero$. By assumption, there exists $\a\in R$ such that
$(\pia,\pia)\not=0$. Since $\rds+S$ and $R$ have the same real
span, we may assume that $\a\in\rds+S$. Also $\pi(\v^{0})\sub
\v^{0}$, so we may assume that $\a=\adot\in\rds$. Set
$$
\b:=\adot+\d_1\in\rds+S\sub\rtimes.
$$
Then $(\pib,\pib)=(\pi(\adot),\pi(\adot))\not=0$. So
$\bt:=\pib\in\pirtimes$. Also
$$
\b+x=\adot+\d_1+\d_1+\d_2\in\adot+2S+S\sub\adot+S\sub\rtimes.
$$
Finally, $\bt+\dt=\pi(\b+x)\in\pirtimes$.

(b) $x\in\rtimes$. By part (i), there exists $\a\in\rtimes$ such
that $(\a,x)\not=0$ and $(\pia,\pia)\not=0$. Now it follows from
the root string property that either $\a+x\in R$ or $-\a+x\in R$.
Take $\bt=\pia$ in the first case and $\bt=\pi(-\a)$ in the
latter case. Then $\bt\in\pirtimes$ and
$\bt+\dt=\pi(\pm\a+x)\in\pirtimes$. This finishes the proof of
lemma.\hfill\qed\vspace{3mm}

Up to now we have only used properties A1--A3 of the automorphism
$\sg$. {\it From now on we assume that $\sg$ satisfies A1--A4}.

As we will see in the next proposition, the conditions A1--A4
imply that $(\gs,\fm\hs)$ satisfies axioms EA1--EA4 of an EALA. It
is shown in [AABGP Chapter I] that most of the structural
properties of an EALA rely on these axioms.

Since $\sg (R)=R$, $\sg$ restricted to $\v$ is an automorphism of
$\v$. Then $\pi(V)$ is the fixed point subspace of $\v$ under
$\sg$. Set
\begin{equation}\label{2-17}
\vs:=\hbox{span}_{\mathbb R}\rs\andd\vpi:=\hbox{span}_{\mathbb
R}\pir= \pi(\v).
\end{equation}
We will see in the next section that in general $\vs$ is a proper
subspace of $\vpi$, so here the upper index $\sg$ does not stand
for the fixed points of $\sg$ on $\v$.

\begin{pro}\label{2-18}
Let $(\gg,\fm,\hh)$ be an EALA and $\sg$ be an automorphism of
$\gg$ such that A1--A4 hold.

(i) If $\rstimes\not=\emptyset$, then $(\gs,\fm,\hs)$ satisfies
axioms EA1--EA4 of an EALA.

(ii) If $\rstimes\not=\emptyset$, then $\rs$ is a reduced EARS in
$\vs$ (in the sense of Definition \ref{EARS}).

(iii) If $\pirtimes\not=\emptyset$, then $\pirtimes$ is an
indecomposable SEARS in $\vpi$.
%and $\pir$ is an irreducible EARS
%in $\vpi$ (in the sense of Definition \ref{EARS}).
\end{pro}

\noindent{\bf Proof.} (i) EA1 holds for $\gs$ by (\ref{2-4}) and
the fact that the form $\fm$ is nondegenerate and invariant on
$\gg$. EA2 holds by A4 and (\ref{2-13}).
%Note that the assumption $\rstimes\not=\emptyset$ implies
%that $\gs\not=\{0\}$. But then condition A4 implies that $\hs\not=\{0\}$,
%since otherwise, $C_{\gs}(\hs)=C_{\gs}(0)=\gs\not=\hs$.
EA3 holds by Lemma \ref{2-14}. Finally EA4 holds by Lemma
\ref{2-15}.

(ii) Since $\rstimes\not=\emptyset$, the form is nontrivial and
positive semidefinite on $\vs$. By (\ref{2-13zz}), R1 holds for
$\rs$. R2 holds by the way $\vs$ is defined. Since $\rs\sub\pir$,
R3 holds by Lemma \ref{2-15}. By part (i), $\gs$ satisfies the
axioms EA1--EA4 of an EALA. Therefore by [AABGP, Theorem 1.29],
the root system $\rs$ of $\gs$ satisfies R4 and R7.

(iii) Clearly the form is positive semidefinite on $\vpi$. By
Lemma \ref{2-16}(iii), $\pir$ and $\pirtimes$ have the same real
span. So it remains to show that if $\at,\bt\in\pirtimes$, then
\begin{equation}\label{2-18a}
2\frac{(\bt,\at)}{(\at,\at)}\in{\bbbz},\andd
\bt-2\frac{(\bt,\at)}{(\at,\at)}\at\in\pir.
\end{equation}
To see (\ref{2-18a}), we use [ABP, Theorem 3.63] (see Remark
\ref{2-19} below), where it is shown that there exists a reduced
irreducible EARS $\rt$ and an isotropic element $\dt$ such that
$$\hbox{span}_{\bbbr}\rt=\vpi\oplus{\bbbr}\dt\andd
\pir\sub\rt\;\;(\mod {\bbbz}\dt).
$$
Now we may consider $\at$ and $\bt$ as elements of $\rt$ $(\mod
{\bbbz}\dt)$ and so (\ref{2-18a}) holds for them as elements of
$\rt$ $(\mod {\bbbz}\dt)$. Thus (\ref{2-18a}) holds for $\at,\bt$
as elements of $\pir$. \hfill\qed\vspace{3mm}

\begin{rem}\label{2-19} {\em The extended affine root
system $\rt$ which appeared in the proof of Proposition \ref{2-17}
is in fact the root system of an EALA $\gt$ which is a covering
algebra of $\gg$ constructed by affinization of the loop algebra
of $\gg$ with respect to $\sg$. In the statement of [ABP, Theorem
3.63], the EALA $(\gg,\fm,\hh)$ is assumed to be tame. However,
what we really need in the proof of Proposition \ref{2-17} is only
the fact that R4 holds for $\rt$. But this holds if axioms
EA1--EA3 of an EALA hold for $\gt$ (see [AABGP, Chapter I]). Now
checking [ABP], one can see that the notion of tameness is not
used in the proof of the fact that EA1--EA3 hold for $\gt$.}
\end{rem}

As in (\ref{riso}), let $\rsiso$ ($\rsniso$) denotes the set of
isolated (nonisolated) isotropic roots of $\rs$.

\begin{cor}\label{2-20} Let $\rstimes\not=\emptyset$. Then
\begin{equation}\label{2-20a}
\rs=(\cup_{i=1}^{k}\rsi)\uplus\rsiso,
\end{equation}
where for each $i$, $\rsi$ is an irreducible reduced EARS in the
real span $\vsi$ of $\rsi$, with $(\rsi,\rs_{j})=\{0\}$ if $i\not=
j$. Moreover if
\begin{equation}\label{2-20-a}
\ri:=\rsi\cup (\la\rsi\ra\cap\rszero),
\end{equation}
then $\ri$ is a reduced indecomposable EARS in $\vsi$.
\end{cor}

\noindent{\bf Proof.} See Proposition \ref{2-18}(ii) and Lemma
\ref{indecom}.\hfill\qed\vspace{3mm}

Note that
\begin{equation}\label{2-10z}
\rsi\sub\ri\andd\hbox{Span}_{\mathbb R}\ri=\hbox{Span}_{\mathbb
R}\rsi=\vsi.
\end{equation}

Since $\rsi$ is an EARS in $\vsi$, we may use the same notation as
in the second paragraph after Definition \ref{EARS}. So $\rsib$ is
an irreducible finite root system in $\vsib$. Fix a choice of a
fundamental system $\piib$ for $\rsib$. Also choose a fixed
preimage $\piid$ in $\rsi$ of $\piib$, under $^{\bar{\;}}$. Let
$\vsid$ be the real span of $\piid$. Then
\begin{equation}\label{2-10w}
\vsi=\vsid\oplus (\vsi)^{0}\andd (\vsi,\vs_{j})=\{0\}\hbox{ for
}i\not=j.
\end{equation}
Then $^{\bar{\;}}$ restricts to an isometry of $\vsid$ onto
$\vsib$. Let
$$
\rsid=\{\at\in\vsid\mid\at+\d\in\rsi\hbox{ for some }\d\in
(\vsi)^{0}\}.
$$
Then $\rsid$ is a finite root system in $\vsid$ which is
isometrically isomorphic under $^{\bar{\;}}$ to $\rsib$. Set
$$
(\rsi)^{0}=\rsi\cap (\vsi)^{0}.
$$
Let $\nu_i=\hbox{dim}(\vsi)^{0}$ and $l_i=\dim\vd^{\sg}_i$. Then
by [AABGP, Chapter II], $(\rsi)^{0}=S_i+S_i$, where $S_i$ is a
semilattice in $(\vsi)^{0}$ of rank $\nu_i$ (see also
\ref{2-15a}). Also the ${\mathbb Z}$--span of $S_i$ is a free
abelian group of finite rank with a basis $B_i\sub S_i$. Now
\begin{equation}\label{2-20aa}
B_i\cup\piid\sub\rsi\sub\hsstar\andd\rsi\sub\la B_i\cup\piid\ra.
\end{equation}
Corresponding to the subspaces $\vsid$ and $(\vsi)^0$ of
$(\hs)^{\star}$ we consider two subspaces $(\vsid)_{\bbbc}$ and
$(\vsi)^{0}_{\bbbc}$ of $\hs$ as follows. Let
\begin{equation}\label{2-20b}
(\vsid)_{\mathbb C}=\sum_{\at\in\piid}{\mathbb
C}t_{\at}\sub\hs\andd (\vsi)^{0}_{\mathbb C}=\sum_{\dt\in
B_i}{\mathbb C}t_{\dt}\sub\hs.
\end{equation}
Set
$$
\vsd_{\bbbc}=\sum_{i=1}^{k}(\vsid)_{\bbbc}\andd
(\vs_{\bbbc})^{0}=\sum_{i=1}^{k}(\vsi)_{\bbbc}^{0}.
$$

\begin{lem}\label{2-100} Supoose $\rstimes\not=\emptyset$.
If $\rsiso=\emptyset$ then
$$
\dot{\v}^{\sg}_{\bbbc}\oplus(\vs_{\bbbc})^{0}=
\sum_{\at\in\rstimes}[\gsat,\gs_{-\at}].
$$
\end{lem}

\noindent{\bf Proof.} Since $\gs$ satisfies axioms EA1-EA4 of an
EALA (see Proposition \ref{2-18}), we have from [AABGP, Chapter
I] that if $\at\in\rstimes$ then
$[\gsat,\gs_{-\at}]={\bbbc}t_{\at}$. So if
$\dot{\a}\in{\piid}\subset\rstimes$, then
$t_{\dot{\a}}\in[\gs_{\dot{\a}},\gs_{-\dot{\a}}]$. If $\d\in S_i$
then from (\ref{2-15a}), we have
$\dot{\a}+\d\in(\rsi)^{\times}\sub\rstimes$ for some
$\dot{\a}\in\dot{R}^{\sg}_{i}$. Then $t_{\dot{\a}+\d}\in
[\gs_{\dot{\a}+\d},\gs_{-\dot{\a}-\d}]$. Thus $t_\d$ is contained
in the ${\bbbc}$-span of $[\gsat,\gs_{-\at}]$, $\at\in\rstimes$
and so $\vsd_{\bbbc}\oplus (\vs_{\bbbc})^{0}$ is contained in the
sum appearing in the statement. Conversely, let $\at\in\rstimes$.
Since $\rsiso=\emptyset$, $\at\in\rsi$ for some $1\leq i\leq k$.
So $\at$ is in the real span of $\piid\cup S_i$. It follows that
$[\gsat,\gs_{-\at}]={\bbbc}t_{\at}\subset\vsd_{\bbbc}\oplus
(\vs_{\bbbc})^{0}$. This completes the
proof.\hfill\qed\vspace{3mm}

Assume that $\rstimes\not=\emptyset$. Since the form $\fm$ is real
valued and positive definite on $\vsid$, it follows that the form
on $\hs$ restricted to $(\vsid)_{\bbbc}$ is nondegenerate, and so
by (\ref{2-10w}) the form on $\vsd_{\bbbc}$ is nondegenerate.
Since the form $\fm$ is nondegenerate on $\hs$, and
$\big(\vsd_{\bbbc}\oplus(\vs_{\bbbc})^{0},(\vs_{\bbbc})^{0}\big)=\{0\}$,
it follows that there exists a subspace $\dd$ of $\hs$ such that
\begin{equation}\label{2-a20}
\begin{array}{l}
\dim{\dd}=\hbox{dim}\sum_{i=1}^{k}(\vs_{i})^{0}_{\mathbb C},\vspace{2mm}\\
\big(\dd,\dd+\vsd_{\bbbc})=\{0\},\vspace{2mm}\\
\fm\hbox{ is nondegenerate on
}\vsd_{\bbbc}\oplus(\vs_{\bbbc})^{0}\oplus\dd.
\end{array}
\end{equation}
Next consider a complement $\w$ of $\vsd_{\bbbc}\oplus
(\vs_{\bbbc})^{0}\oplus\dd$ in $\hs$ such that
\begin{equation}\label{2-b20}
\begin{array}{l}
(\vsd_{\bbbc}\oplus (\vs_{\bbbc})^{0}\oplus\dd,\w)=\{0\},\vspace{2mm}\\
\fm\hbox{ is nondegenerate on }\w.
\end{array}
\end{equation}
For $1\leq j\leq k$, consider a subspace $\dj$ of $\dd$ such that
\begin{equation}\label{2-20d}
\begin{array}{l}
\hbox{dim}(\vs_{j})^{0}_{\mathbb C}=\hbox{dim}\dj, \vspace{2mm}\\
\fm\hbox{ is nondegenerate on
}(\dot{\v}^{\sg}_{j})_{\bbbc}\oplus(\vs_{j})^{0}_{\bbbc}\oplus\dj.
\end{array}
\end{equation}

Set
\begin{equation}\label{2-20z}
\hsi=(\vsid)_{\mathbb C}\oplus (\vsi)^{0}_{\mathbb C}\oplus\di.
\end{equation}
Then the form on $\hsi$ is nondegenerate and
$$
\hs=\sum_{i=1}^{k}\hsi\oplus\w.
$$

%Let $\mi$ be the subalgebra of $\gs$ generated by root spaces
%$\gs_{\at}$ with $\at\in R_{i}\setminus\{0\}$.
Put
$$
\gsi=\hsi\oplus\sum_{\at\in R_{i}\setminus\{0\}}\gsat
$$
\begin{pro}\label{2-23}
(i) $\gsi$ is a subalgebra of $\gs$.

(ii) $\hsi$ is an abelian subalgebra of $\gsi$. Moreover the form
restricted to $\hsi$ is nondegenerate.

(iii) $C_{\gsi}(\hsi)=\hsi$.
\end{pro}

\noindent{\bf Proof.} (i) Set
$$\ki=\sum_{\at\in\ri\setminus\{0\}}\gsat\andd
\ti=\sum_{\at\in\ri}{\mathbb C}t_{\at}=\sum_{\at\in\rsi}{\mathbb
C}t_{\at}=\sum_{\at\in (\rsi)^{\times}}{\mathbb C}t_{\at},
$$
and let $\mi$ be the subalgebra of $\gs$ generated by $\ki$. By
(\ref{2-20aa}) and (\ref{2-20b}),
$$
\ti=(\vsid)_{\mathbb C}\oplus (\vsi)^{0}_{\mathbb C}.
$$
We first claim that $\mi=\ki\oplus\ti$. To see this note that by
Proposition \ref{2-18}, EA1--EA4 hold for $\gs$, therefore for any
$\at\in\rs$,
$$
[\gsat,\gs_{-\at}]={\mathbb C}t_{\at}
$$
(See [AABGP, Chapter I]). So $\ti\sub\mi$. Thus
$\ti\oplus\ki\sub\mi$. To see the equality it is enough to show
that $\ti\oplus\ki$ is closed under $[\;,\;]$. But $\ti\sub\hs$,
so $[\ti,\ti\oplus\ki]\sub\ki$. To see
$[\ki,\ki]\sub\ti\oplus\ki$, let $\at,\bt\in\ri\setminus\{0\}$,
and $[\gsat,\gs_{\bt}]\not=\{0\}$. In particular $\at+\bt\in\rs$.
If $\bt=-\at$, then $[\gsat,\gs_{\bt}]={\mathbb
C}t_{\at}\sub\ti$. So we may assume that $\bt+\at\not=0$. If
$\at+\bt$ is nonisotropic, then either $(\at+\bt,\at)\not=0$ or
$(\at+\bt,\bt)\not=0$. In either cases $\at+\bt$ is a
nonisotropic root in $\rs$ which is not orthogonal to
$\at\in\rsi$ or $\bt\in\rsi$, so $\at+\bt\in (\rsi)^{\times}$.
Thus
$$
[\gsat,\gs_{\bt}]\sub\gs_{\at+\bt}\sub\ki.
$$
If $\at+\bt$ is isotropic, then $\at+\bt\in
\la\rsi\ra\cap\rszero\sub\ri$. Thus
$[\gsat,\gs_{\bt}]\sub\gs_{\at+\bt}\sub\ki$, as
$\at+\bt\in\ri\setminus\{0\}$. This completes the proof of our
claim.

Next note that
$$
\gsi=\hsi\oplus\sum_{\at\in\ri\setminus\{0\}}\gsat
=\ti\oplus\ki\oplus\di=\mi\oplus\di.
$$
Since $\di\sub\hsi\sub\hs$, $[\di,\di\oplus\ki]\sub\ki$. But
$\mi$ is the subalgebra generated by $\ki$, so
$[\di\oplus\mi,\di\oplus\mi]\sub\mi$.

(ii) The first assertion is clear. It follows from (\ref{2-20d})
that the form is nondegenerate.

(iii) Suppose to the contrary that $C_{\gsi}(\hsi)\not\sub\hsi$.
Then there exists $x=x_{0}+x_{\at_1}+\cdots +x_{\at_t}\in
C_{\gsi}(\hsi)$ such that $x_0\in\hsi$, $x_{\at_j}$'s are nonzero
and $\at_j$'s are distinct roots of $\ri\setminus\{0\}$. Then for
any $h\in\hsi$,
$$
0=[h,x]=\at_{1}(h)x_{\at_{1}}+\cdots +\at_{t}(h)x_{\at_{t}}.
$$
Thus $\at_{j}(h)=0$ for $1\leq j\leq t$, $h\in\hsi$. That is
$(t_{\at_{j}},\hsi)=\{0\}$. Since $t_{\at_{j}}\in\ti\sub\hsi$, it
follows from part (ii) that $\at_{j}=0$ for all $j$, which is a
contradiction.\hfill\qed\vspace{3mm}

\begin{pro}\label{2-24}
Suppose that $\rstimes\not=\emptyset$. Then $(\gsi,\fm,\hsi)$
satisfies axioms EA1--EA5(a) of an EALA.
\end{pro}

\noindent{\bf Proof.} By Proposition \ref{2-18}, EA1--EA2 hold
for $\gs$. So for any $\at\in\rs$, the form restricted to
$\gs_{\at}\oplus\gs_{-\at}$ is nondegenerate. In particular, this
holds for any $\at\in\ri$. Thus $\fm$ is nondegenerate on $\gsi$
and EA1 holds for $\gsi$. From (\ref{2-20aa}) we have
$\ri\subset(\hs)^{\star}$. Now from the way $\gsi$ is defined it
follows that $\hbox{ad}_{\gsi}h$ is diagonalizable for all
$h\in\hsi$. This together with Proposition \ref{2-23}(ii) imply
that EA2 holds for $\gsi$. EA3 holds for $\gsi$ as it holds for
$\gs$. EA4 holds for the root system $\rsi$ of $\gsi$ since
$\rsi\sub\rs$ and EA4 holds for $\rs$. Finally EA5(a) holds by
Corollary \ref{2-20}.\hfill\qed\vspace{3mm}

Set
\begin{equation}\label{2-24z}
\i=\sum\gs_{\dt},
\end{equation}
where sum runs through
\begin{equation}\label{2-24w}
\dt\in\rsiso\setminus(\cup_{i=1}^{k}\ri)=
\rsiso\setminus\big(\cup_{i=1}^{k}(\ri)_{\hbox{iso}}\big).
\end{equation}
%Let
%\begin{equation}\label{2-25}
%\z=\w\oplus\i.
%\end{equation}
Then
\begin{equation}\label{2-14z}
\gs=\hs\oplus\sum_{\at\in\rs\setminus
\{0\}}\gsat=\sum_{i=1}^{k}\gsi\oplus\w\oplus\i.
\end{equation}

For our next proposition we need to state some results regarding
the core $\gsc$ of $\gs$. In particular, we would like to obtain
some criteria for the tameness of $\gs$. The notions of core and
tameness are defined for an EALA but as we will see in next
section $\gs$ may not be in general an EALA. In the sequel, we
also have some other Lie algebras which satisfy EA1--EA4 but are
not in general an EALA. So let us define the {\it core} $\gc$ of a
triple $(\gg,\fm,\hh)$ satisfying EA1--EA4 to be the subalgebra
of $\gg$ generated by root spaces $\gg_{\a}$, $\a\in\rtimes$.
$\gg$ is called {\it tame} if the centralizer of $\gc$ in $\gg$ is
contained in $\gc$. It follows that $\gc$ is a perfect ideal of
$\gg$.

Since $\sg(\gg_{\a})=\gg_{\sg(\a)}$, we have $\sg(\gc)=\gc$. Thus
\begin{equation}\label{2-16a}
\gc=\bigoplus_{i=0}^{m-1}(\gc)_{i},
\end{equation}
where $(\gc)_{i}=\gc\cap\gg_{\bar{i}}$.

\begin{lem}\label{2-16b} Suppose $\rstimes\not=\emptyset$. Then

(i) $\gc\cap\gs$ is the sum of spaces of the forms:
\begin{equation}\label{2-16c}
\gspia,\quad \a\in R,\quad \pia\in\pirtimes,
\end{equation}
and
\begin{equation}\label{2-16d}
[\gipia,\gg_{-\bar{i},\pib}],\;\;\a,\b\in R,\;\;
\pia,\pib\in\pirtimes,\;\;i\in{\bbbz}.
\end{equation}

(ii) $\gsc$ is the sum of spaces of the forms (\ref{2-16c}) and
their commutators.

(iii) $\gsc\sub\gs\cap\gc$.
%Moreover, suppose that if $\a,\b\in
%R$ with $\pia,\pib\in\pirtimes$ and $\pia+\pib\in\pir^{0}$, then
%\begin{equation}\label{2-16f}
%[\gipia,\gg_{-\bar{i},\pib}]\not=\{0\}\Longrightarrow\gs_{\pi(\a+\b)}\sub\gsc.
%\end{equation}
%Then $\gsc=\gs\cap\gc$.

(iv) If $\gsc\cap\hs=\gc\cap\hs$, then
$C_{\gs}(\gs\cap\gc)\subseteq\gs\cap\gc^{\perp}$. In particular if
$\gg$ is tame and $\gsc=\gs\cap\gc$, then $\gs$ is tame.

%(v) If $\gsc\cap\hs\subsetneqq\gc\cap\hs$, then $\gs$ is not tame.
\end{lem}

\noindent{\bf Proof.} (i)This is an immediate consequence of
(\ref{2-16a}) and [ABP, Lemma 3.53].

(ii) Let $S$ be the sum of the spaces in (\ref{2-16c}) and their
commutators. Then $S\sub\gsc$. To show the reverse inclusion it
is enough to show that $S$ is closed under bracket, since $S$
contains all the generators of $\gsc$. So it is enough to show
that for any three spaces $\gspia$, $\gspib$, $\gs_{\pi(\gamma)}$
as in (\ref{2-16c}),
\begin{equation}\label{2-16e}
[\gspia,[\gspib,\gspig]]\sub S.
\end{equation}
Certainly, we may assume that the two brackets involved in
(\ref{2-16e}) are nonzero. Then by (\ref{2-13a}),
$\pib+\pig=\pi(\a')$ for some $\a'\in R$ and $[\gspib,\gspig]\sub
\gs_{\pi(\a')}$. Now if $\pi(\a')\in\pirtimes$, then
$$
[\gspia,[\gspib,\gspig]]\sub[\gspia,\gs_{\pi(\a')}]\sub S.
$$
If $\pi(\a')\in\pir^{0}$, then $\pia+\pi(\a')$ is nonisotropic and
by (\ref{2-13a}), $\pia+\pi(\a')=\pi(\b')$ for some $\b'\in R$,
and $[\gspia,\gs_{\pi(\a')}]\sub\gs_{\pi(\b')}$. But
$\gs_{\pi(\b')}\sub S$, since $\pi(\b')\in\pirtimes$. This
completes the proof of part (ii). Now (iii) is an immediate
consequence of (i) and (ii)

(iv) Let $x\in C_{\gs}(\gs\cap\gc)$. Then
$x=\sum_{i=0}^{n}x_{\ati}$ where $\ati$ are distinct roots of
$\rs$ with $\at_0=0$, and $x_{\ati}\in\gs_{\ati}$ for $i\not=0$.
Let $\at,\bt\in\rstimes$. By part (i),
\begin{equation}\label{2-16g}
\sum_{i=0}^{n}[x_{\ati},\gs_{\at}]=[x,\gs_{\at}]=\{0\}=\sum_{i=0}^{n}
[x_{\ati},[\gg_{\bar{j},\at},\gg_{-\bar{j},\bt}]]=
[x,[\gg_{\bar{j},\at},\gg_{-\bar{j},\bt}]],
\end{equation}
for $j\in{\bbbz}$, where
$[x_{\ati},\gs_{\at}]\subseteq\gs_{\ati+\at}$ and
$[x_{\ati},[\gg_{\bar{j},\at},\gg_{-\bar{j},\bt}]]\subseteq\gs_{\ati+\at+\bt}$.
It follows that for $0\leq i\leq n$,
$$
[x_{\ati},\gs_{\at}]=\{0\}\andd[x_{\ati}[\gs_{\bar{j},\at},\gs_{-\bar{j},\bt}]]=\{0\},
$$
and so $x_{\ati}\in C_{\gs}(\gs\cap\gc)$ for $0\leq i\leq n$.
Therefore we must show that for each $i$,
$x_{\ati}\in\gc^{\perp}$. By (\ref{2-3}), it is enough to show
that for each $i$, $(x_{\ati},\gs\cap\gc)=\{0\}$. Since the form
is nondegenerate and invariant on $\gs$, and since $\gs_{0}=\hs$,
it follows that
\begin{eqnarray*}
(x_{0},\gs\cap\gc)&=&(x_{0},\hs\cap\gc)=(x_0,\hs\cap\gsc)\\
&=& (x_{0},\gsc)=(x_{0},[\gsc,\gsc])\\
&=&([x_0,\gsc],\gsc)\subseteq([x_0,\gs\cap\gc],\gsc)=\{0\}.
\end{eqnarray*}
Thus $x_{0}\in \gc^{\perp}$. Next, we show that for $1\leq i\leq
n$, $x_{\ati}\in\gc^{\perp}$ (note that $\ati\not=0$). Fix $1\leq
i\leq n$. Let $\at,\bt\in\rstimes$ and $\gamt=\at+\bt$. Let
$y_{\gamt}\in\gs_{\gamt}$ or $y_{\gamt}\in
[\gg_{\bar{j},\at},\gg_{-\bar{j},\bt}]$ for some $j\in{\bbbz}$.
We must show that $(x_{\ati},y_{\gamt})=0$. If $\ati+\gamt\not=0$
then by (\ref{2-13z}), $(x_{\ati},y_{\gamt})=0$. So we may assume
that $\gamt=-\ati$. Since $y_{\gamt}\in\gs\cap\gc$, and since
EA1-EA2 hold for $\gs$, we see from [AABGP, I.(1,10)] that
$$
0=[x_{\ati},y_{\gamt}]=(x_{\ati},y_{\gamt})t_{\ati}.
$$
Thus $(x_{\ati},y_{\gamt})=0$. This completes the proof of the
first statement in (iv). Now if $\gg$ is tame then
$\gc^{\perp}\subseteq\gc$. So the second statement of (iv) is an
immediate consequence of the first statement.
%(v) From the structure of an EARS we know that each element of

%$R$ if of the form $\dot{\a}+\d$ where $\dot{\a}$ belongs to some

%finite root system $\rd$ and $\d$ belongs to a lattice of finite

%rank with some basis $\d_1,\ldots,\d_{\nu}\in R^{0}$. Moreover,

%$R$ contains a basis $\dot{\a}_1,\ldots,\dot{\a}_{\ell}$ of

%$\rd$. Take

%$\vd_{\bbbc}=\sum_{i=1}^{\ell}{\bbbc}t_{\dot{\a}_{i}}$ and

%$\v^{0}_{\bbbc}=\sum_{i=1}^{\nu}{\bbbc}t_{\d_{i}}$. It follows

%that $\hh\cap\gc=\vd_{\bbbc}\oplus\v^{0}_{\bbbc}$ (see [A3] or

%[BGK]). In particular $(h,\vd_{\bbbc}\oplus\v^{0}_{\bbbc})=\{0\}$

%for all $h\in\v^{0}_{\bbbc}$, and $(h,h)\not=0$ for all nonzero

%elements $h\in\vd_{\bbbc}$. Now it follows from the assumption in

%the statement that there exists

%$$0\not=h\in (\gc\cap\hs)\setminus(\gsc\cap\hs)\quad\hbox{with}

%\quad(h,\gsc\cap\hs)=\{0\}.

%$$

%Now if $\pia\in\rstimes$ and $x\in\gspia$, then

%$t_{\pia}\in\gsc\cap\hs$ and so we have

%$$

%[h,x]=\pia(h)x=(t_{\pia},h)=0.

%$$

%Thus $h\in C_{\gs}(\gsc)\setminus\gsc$ which implies that $\gs$ is

%not tame. This completes the proof.
\hfill\qed\vspace{3mm}

\begin{rem}\label{2-16-g}
{\em As we will see in the next section, in many examples $\gsc$
is a proper subspace of $\gs\cap\gc$.}
\end{rem}

Note that from Proposition \ref{2-16b}(ii), we have
\begin{equation}\label{2-25a}
\gsc\sub\sum_{i=1}^{k}\gsi\sub\gs.
\end{equation}

\begin{pro}\label{2-26} (i) $\w\sub C_{\gs}(\sum_{i=1}^{k}\gsi)$.

(ii) $\sum_{\dt\in\rsiso}\gs_{\dt}\sub C_{\gs}(\gsc)$. In
particular $\w\oplus\i\sub C_{\gs}(\gsc)$.

%(iii) For each $i$, $\i\sub C_{\gs}(\gsi/\di)$.
\end{pro}

\noindent{\bf Proof.} (i) Since $\w\sub\hs$, we have
$[\w,\hs]=\{0\}$. Therefore, it remains to show that
$[\w,\gsat]=\{0\}$ for any $\at\in\ri\setminus\{0\}$, $1\leq
i\leq k$. Let $\at\in\ri\setminus\{0\}$. Then $\at\in\la\rsi\ra$.
By Corollary \ref{2-20} $\rsi$ is tame, so $\la\rsi\ra=\la
(\rsi)^{\times}\ra$. By (\ref{2-20aa}), (\ref{2-20b}) and
(\ref{2-b20}), we have $\w\sub\hbox{ker}\bt$, for any $\bt\in
(\rsi)^{\times}$. Therefore $[\w,\gsat]=\{0\}$.

(ii) Let $\dt\in\rsiso$. Since $\dt+\at\not\in\rs$ for any
$\at\in\rstimes$, $[\gsat,\gs_{\dt}]=\{0\}$. But $\gsc$ is
generated by root spaces $\gsat$, $\at\in\rstimes$, so
$[\gs_{\dt},\gsc]=\{0\}$. Now the proof is complete, using part
(i) and (\ref{2-25a}). \hfill\qed\vspace{3mm}

\begin{cor}\label{2-a27}
If $\rstimes\not=\emptyset$, then for each $i$, the triple
$$(\gsi\oplus\w,\fm,\hsi\oplus\w)$$ satisfies EA1-EA5(a).
\end{cor}

\noindent{\bf Proof.} By Proposition \ref{2-26}, $\gsi\oplus\w$
is a Lie algebra. By Proposition \ref{2-24} the form $\fm$ is
nondegenerate on $\gsi$, so by (\ref{2-b20}) and (\ref{2-20d}),
the form $\fm$ is nondegenerate on the Lie algebra
$\gsi\oplus\w$. By Propositions \ref{2-23}(iii) and
\ref{2-26}(ii), $\hsi\oplus\w$ is self-centralizing in
$\gsi\oplus\w$. By acting as zero on $\w$ we may identify
elements of $(\hsi)^{\star}$ as elements of
$(\hsi\oplus\w)^{\star}$. Then it follows from Proposition
\ref{2-24} that EA2 holds for $\gsi\oplus\w$ and that $\gsi$ and
$\gsi\oplus\w$ have the same root system. Now it is clear that
EA3-EA5(a) hold for $\gsi\oplus\w$.\hfill\qed\vspace{3mm}

Recall that the EALA $\gg$ is called nondegenerate if the real
space $\v^0$ and its complex span in $\hh^\star$ have the same
dimension.
%$\dim\v^{0}=\dim(\v^{0}\otimes_{\bbbr}{\bbbc})$.
We summarize our results in the following theorem.

\begin{thm}\label{2-27} Let $(\gg,\fm,\hh)$ be an EALA with
corresponding root system $R$. Let $\sg$ be an automorphism of
$\gg$ satisfying A1--A4. Let $\gs$ ($\hs$) be the fixed point
subalgebra of $\gg$ ($\hh$), under $\sg$, and let $\rs$ be the
root system of $\gs$ with respect to $\hs$. Suppose
$\rstimes\not=\emptyset$, then

(i) $(\gs,\fm,\hs)$ satisfies axiom EA1--EA4 of an EALA and $\rs$
is a reduced EARS.

(ii) $\rs=\big(\cup_{i=1}^{k}\ri\big)\cup\rsiso$ where for each
$i$, $\ri$ is an indecomposable reduced EARS with
$(\ri^{\times},R_{j}^{\times})=\{0\}$ if $i\not=j$. (The union is
not necessarily disjoint.)

(iii) $\hs=\sum_{i=1}^{k}\hsi\oplus\w$ where $\hsi$ and $\w$ are
some subspaces of $\hs$ with $(\hsi,\w)=\{0\}$ and
$[\w,\gsi]=\{0\}$, for each $i$.

(iv) $\gs=\sum_{i=1}^{k}\gsi\oplus\w\oplus\i$, where for each $i$,
$(\gsi,\fm,\hsi)$ is a Lie algebra satisfying axioms EA1--EA5(a)
of an EALA, and $\i$ is a subspace of $\gs$, satisfying $\i\sub
C_{\gs}(\gsc)$. Moreover, $\i=\{0\}$ if $\riso=\emptyset$.

(v) If $i\not=j$, then $[(\gs_{i})_c,(\gs_{j})_c]=\{0\}$. Moreover
$\gsc=\sum_{i=1}^{k}(\gs_i)_c$. ($(\gs_{i})_c$ is the core of
$\gsi$.)

(vi) If $\gg$ is nondegenerate then for each $i$,
$\dim\hsi=l_i+2\nu_{i}$.
\end{thm}

\noindent{\bf Proof.} For (i) see Proposition \ref{2-18}. For
(ii) see Corollary \ref{2-20}. For (iii) see (\ref{2-b20}) and
Proposition \ref{2-26}. For (iv) see (\ref{2-24z}), Proposition
\ref{2-24}, and Proposition \ref{2-26}.

For the first statement of (v), it is enough to show that if
$\pia\in R_{i}^{\times}$ and $\pib\in R_{j}^{\times}$, then
$\pia+\pib$ is not a root in $\rs$. By part (ii),
$(\pia,\pib)=0$. Therefore $\pia+\pib$ is not orthogonal to both
$\pia$ and $\pib$ and so cannot be a root (by part (ii)). This in
particular shows that $\sum_{i=1}^{k}(\gs_{i})_{c}$ is a Lie
algebra which contains all generators of $\gsc$. So the second
statement of (v) holds. (vi) is clear from the construction of
$\hsi$. \hfill\qed\vspace{3mm}

As a corollary we can state a weak version of a result which is
due to [BM] (see Remark \ref{2-28-1}).
\begin{cor}\label{2-29}
Let $(\gg,\fm,\hh)$ be a finite dimensional complex simple Lie
algebra, where $\hh$ is a Cartan subalgebra and $\fm$ is the
Killing form on $\gg$. Let $\sg$ be an automorphism of $\gg$
satisfying A1-A2 and A4. Then $(\gs,\fm,\hs)$ is a reductive Lie
algebra.
\end{cor}

\noindent{\bf Proof.} It is easy to see that A3 is a consequence
of A1. $\gg$ is an EALA of nullity zero and so there is no
nonzero isotropic root. Therefore by (\ref{2-24z}) and
(\ref{2-24w}), $\i=\{0\}$. By Corollary \ref{2-20} and part (ii)
of Theorem \ref{2-27}, $\ri=\rsi$ is an irreducible finite roots
system in $\vsi$. By Theorem \ref{2-27}(iv) and (\ref{2-20z}),
each $\gsi$ is an EALA of nullity zero, where $\hsi$ is of
dimension equal to the rank of $\rsi$. Thus $\gsi$ is a finite
dimensional simple Lie algebra. Finally $\w$ is an abelian
subalgebra of $\gg$. In fact it follows from Proposition
\ref{2-26}(i) that $\w={\mathcal Z} (\gs)$.\hfill\qed\vspace{3mm}

\begin{rem}\label{2-28-1}{\em (i) According to [BM], if $\sg$
satisfies only A1, then $\gs$ is a reductive Lie algebra. A
version of conjugacy for Cartan subalgebras is used in the proof
and no such result is known for a general EALA.

(ii) According to Theorem \ref{2-27} and Corollary \ref{2-a27}, we
may express the result of [BM] in another way. Namely, if $\sg$
is a finite order automorphism of a finite dimensional complex
simple Lie algebra $\gg$, then $\gs={\mathcal J}_1\oplus{\mathcal
J}_2$ where ${\mathcal J}_1$ is a semisimple Lie algebra and
${\mathcal J}_{2}$ is an EALA such that ${\mathcal
J}_{2}/{\mathcal Z}({\mathcal J}_{2})$ is a finite dimensional
simple Lie algebra.

(iii) Similar to part (ii), we may state Theorem \ref{2-27}(iv)
in a different way. In fact by Corollary \ref{2-a27}, for a fixed
$i$ we can consider $\gsi\oplus\w$ as a Lie algebra satisfying
EA1-EA5(a). So part (iv) of the Theorem can be restated to say
$\gs$ is a direct sum of some EALAs and a subspace $\i$
satisfying $\i\sub C_{\gs}(\gsc)$.

(iv) One may define a notion of nondegeneracy for a triple
$(\gg,\fm,\hh)$ satisfying EA1-EA5(a) exactly in the same way
which one defines this for an EALA. Therefore by Theorem
\ref{2-27} for each $i$, $\gsi$ is nondegenerate. }
\end{rem}

\begin{cor}\label{2-28}
Let $(\gg,\fm,\hh)$, $\sg$ and $\rs$ be as in Theorem \ref{2-27}.
If $\rs$ is an irreducible EARS, then $(\gs,\fm,\hs)$ is an EALA.
\end{cor}

\noindent{\bf Proof.} By Theorem \ref{2-27}, $\gs$ satisfies
EA1--EA5(a). Now EA5(b) also holds as $\rs$ is
irreducible.\qed\vspace{3mm}

\begin{cor}\label{2-30} Let $(\gg,\fm,\hh)$ be as in Theorem
\ref{2-27} and assume that $\rstimes\not=\emptyset$ and that $\gg$
is tame. If $\gs\cap\gc=\gsc$, then $\gs=\sum_{i=1}^{k}\gsi$,
where each $\gsi$ is an EALA.
\end{cor}

\noindent{\bf Proof.} (i) By Lemma \ref{2-16b}(iv),
$C_{\gs}(\gsc)\sub\gsc$. Therefore by Theorem \ref{2-27},
$(\gs,\fm,\hs)$ is a tame Lie algebra satisfying EA1--EA4. Now it
follows from [ABP, Lemma 3.62] that $\gs$ also satisfies EA5(b)
(that is $\rsiso=\emptyset$). In particular $\i=\{0\}$ and each
$\ri$ is an irreducible reduced EARS. Therefore, by part (iv) of
Theorem \ref{2-27} each $\gsi$ is an EALA. Since $\gs$ is tame we
have from Proposition \ref{2-26} that $\w\subset\gsc\cap\hs$. By
Lemma \ref{2-100},
$$\gsc\cap\hs=\sum_{\at\in\rstimes}[\gsat,\gs_{-\at}]=
\vsd_{\bbbc}\oplus(\vs)^{0}_{\bbbc}.
$$
But $\vsd_{\bbbc}\oplus(\vs)^{0}_{\bbbc}$ and $\w$ have zero
intersection. Thus $\w=\{0\}$. \hfill\qed\vspace{4mm}
\vspace{.5cm}

\section{EXAMPLES}
In this section we present several examples which elaborate on the
results obtained in Section 2. In \ref{3-2}-\ref{3-5} below a
large class of examples is introduced which illustrate how the
terms $\gsi$'s, $\w$ and $\i$ (see Theorem \ref{2-27}) appear as
the fixed points of automorphisms. In \ref{3-6}-\ref{3-9} and
\ref{3-11} we recall some examples from [A2], in which certain
finite order automorphisms of EALA are given. We show that these
automorphisms satisfy conditions A1-A4. Then, using Theorem
\ref{2-27}, we are able to give a new proof of the results
obtained in [A2], namely to prove that many examples of EALA (of
types $D_{\ell}$, $A_1$, $B_{\ell}$, $C_{\ell}$, and $BC_{\ell}$)
can be obtained as the fixed points of automorphisms of some
other EALA (of types $A_{\ell}$, $D_{\ell}$ and $C_{\ell}$) which
may have a simpler structure. Finally in \ref{3-12}, we present
an example regarding the results in [ABP].

\begin{exa}\label{3-2}
Let $(\gg,\fm,\hh)$ be an EALA of type $X_{\ell}$ and nullity
$\nu$. Let $R$ be the corresponding irreducible extended affine
root system and denote its root lattice by $Q$. Then
$$Q={\mathbb
Z}\a_{1}\oplus\cdots\oplus{\bbbz}\a_{\ell}\oplus{\mathbb
Z}\d_{1}\oplus\cdots\oplus{\bbbz}\d_{\nu},
$$
where $\{\a_{1},\ldots,\a_{\ell}\}$ is a subset of $R$ which form
a basis for a finite root system of type $X_{\ell}$ and
$\{\d_{1},\ldots,\d_\nu\}$ is a subset of $R$ which form a basis
of the radical of the form restricted to the real span of $R$.
Consider any group homomorphism $\phi:Q\longrightarrow{\mathbb
C}\setminus\{0\}$. $\phi$ is uniquely determined by specifying
$\phi(\a_{i})$ for $1\leq i\leq\ell$ and $\phi(\d_j)$ for $1\leq
j\leq \nu$. The homomorphism $\phi$ induces an automorphism $\sg$
of $\gg$ by letting
$$
\sg_{|_{\gg_a}}=\phi(\a)\hbox{id}_{\gg_a}\quad\hbox{for }\a\in R.
$$
Since $\gg_{0}=\hh$ and $\phi(0)=1$, we have $\sg(h)=h$ for all
$h\in\hh$.  Note that $\sg$ is of finite order if and only if
$\phi(\a_i)$'s for $1\leq i\leq\ell$ and $\phi(\d_{j})$'s for
$1\leq j\leq\nu$ are roots of unity. Assume that $\sg^{m}=
\hbox{id}$, for some positive integer $m$. Let $\gs$ be the set of
fixed points of $\sg$. Clearly $\sg(\hh)=\hh$ and $\hs=\hh$. Since
$$(\gg_{\a},\gg_{\b})=\{0\}\quad\hbox{unless}\quad\a+\b=0,
$$
it follows that
$$\big(\sg(x),\sg(y)\big)=(x,y)\quad\hbox{for all }x,y\in\gg.
$$
Since $\sg(\a)=\a$ for all $\a\in(\hs)^\star$, it follows from
[ABP, Proposition 3.25] that $C_{\gs}(\hh)=\hh$. Therefore $\sg$
satisfies conditions A1--A4. Thus by Theorem \ref{2-27}, $\gs$
satisfies axioms EA1--EA4. Therefore if $\rstimes\not=\emptyset$,
then
\begin{equation}\label{3-2e}
\gs=\sum_{i=1}^{k}\gsi\oplus\w\oplus\i,
\end{equation}
where for each $i$, $(\gsi,\fm,\hsi)$ satisfies EA1--EA5(a) and
$\w$ and $\i$ are as in Theorem \ref{2-27}. Note that
$\gs\cap\gg_{\a}\not=\{0\}$ if and only if $\phi(\a)=1$. In
particular,
$$
\rs=\{\a\in R\mid\gs\cap\gg_\a\not=\{0\}\}=\{\a\in
R\mid\gg_{\a}\sub\gs\}=\{\a\in R\mid\phi(\a)=1\}.
$$
\end{exa}

%\begin{exa}\label{3-3} In the previous example suppose that $\gg$
%is a finite dimensional complex simple Lie algebra. Then by
%Corollary \ref{2-29}, $\gs$ is a reductive Lie algebra.
%\end{exa}

\begin{exa}\label{3-4}
In Example \ref{3-2} assume that $X_{\ell}$ has one of the types
$A_{\ell}$, $D$, $E$, $B_{\ell}$, $C_{\ell}$, $F_4$ or $G_2$. Let
$S$ and $L$ be semilattices which appear in the structure of $R$
(see (\ref{2-15a})). We want to impose some restrictions on the
semilattices $S$ and $L$ as follows. For types $A_{\ell}$
$(\ell\geq 2)$, $D$ and $E$ we impose no restriction as we know
from [AABGP, II.2.32] that for these types the semilattice $S$ is
always a lattice. For type $A_{1}$ assume that $S$ is a lattice
and for the remaining types assume that both $S$ and $L$ are
lattices and that $S=L$. (The root systems of toroidal Lie
algebras, with some derivations added, are of this form). We
claim that for such types and under the above restrictions the
axiom EA5(b) also holds

Suppose that $\rstimes\not=\emptyset$. So there exists $\a\in
R^{\times}$ such that $\phi(\a)=1$. Let $\d\in\rszero$. That is
$\phi(\d)=1$. Under the above assumptions, it follows from the
structure of EARS of type $X_{\ell}$ that $\a+\d\in R$. Now
$\phi(\a+\d)=1$ and so $\a+\d\in\rs$. This shows that EA5(b)
holds. Thus by Example \ref{3-2},
\begin{equation}\label{3-4e}
\gs=\sum_{i=1}^{k}\gsi\oplus\w,
\end{equation}
where each $\gsi$ is an EALA and $\w$ is an abelian subalgebra of
$\gg$.
\end{exa}

\begin{exa}\label{3-5}
Let $\gg$ and $\sg$ be as in Example \ref{3-2} and let
$X_{\ell}=A_1$. We have
$$
R=(S+S)\cup(\pm\dot{\a}+S),
$$
where $S$ is a semilattice in the real span $\v^{0}$ of $R^{0}$.
If $S$ is a lattice, then by Example \ref{3-4} axiom EA5(b) holds
for $\gs$.

Next suppose that $S$ is not a lattice. We show that in this case
it might happen that the axiom EA5(b) does not hold. To see this
let the nullity $\nu$ of $R$ be $3$. Let
$$
S=\{\sum_{i=1}^{3}m_i\d_i\mid m_i\in{\mathbb Z}\hbox{ and
}m_{i}m_{j}\equiv 0\hbox{ mod }2,\hbox{ if }i\not=j\}.
$$
Then $S$ is a semilattice in $\v^{0}$ and $Q={\mathbb
Z}\dot{\a}\oplus{\mathbb Z}\d_{1}\oplus{\mathbb
Z}\d_{2}\oplus{\mathbb Z}\d_{3}$. Define
$$
\phi(\dot{\a})=-1,\quad\phi(\d_1)=1,\quad\phi(\d_{2})=1,\andd\phi(\d_{3})=-1.
$$
Now $\phi(\d_1+\d_2)=1$ and so $\d:=\d_1+\d_2\in\rszero$. We
claim that $\d$ is isolated, that is $\d\in\rsiso$. Suppose to
the contrary that there exists
$\a=\pm\dot{\a}+m_1\d_1+m_2\d_2+m_3\d_3\in\rs$ such that
$\a+\d\in\rs$. We have
$$
1=\phi(\a)=-\phi(\d_3)^{m_3}=-(-1)^{m_3},
$$
and so $m_3=2k_3+1$ for some $k_3\in{\mathbb Z}$. Since
$\a\in\rs\sub R$, we must have $m_1=2k_1$ and $m_2=2k_2$ for some
$k_1,k_2\in{\mathbb Z}$. Then
$$
\a+\d=\pm\dot{\a}+(2k_1+1)\d_1+(2k_2+1)\d_2+(2k_3+1)\d_3\in\rs\sub
R.
$$
But this contradicts the fact that $R$ contains no such root.
\end{exa}

For our next few examples we need the following setting (see
[BGK], [AABGP] and [A2]). Let $\nu\geq 1$. Let
$\be=(e_{1},\ldots,e_\nu)$ be a vector in ${\mathbb C}^{\nu}$ and
let $\bq=(q_{ij})$ be a $\nu\times\nu$--matrix such that
\begin{equation}\label{setting1}
e_i=\pm 1,\;\;\;q_{ii}=1\hbox{ for }1\leq i\leq\nu,\andd
q_{ij}=q_{ji}\hbox{ for }1\leq i\not= j\leq\nu.
\end{equation}
Let $\aa$ be the associative algebra over ${\mathbb C}$ with
generators $x_{i}, x_{i}^{-1}$ subject to the relations
$x_{i}x_{j}=q_{ij}x_{j}x_{i}$. Then
$$\aa=\bigoplus_{\d\in\zn}{\bbbc}\xd,\hbox{ where }
\xd=x_{1}^{n_{1}}\cdots x_{\nu}^{n_{\nu}}\hbox{ for }
\d=(n_1,\ldots,n_\nu)\in\zn.
$$
Let $\bar{\;}$ be the involution on $\aa$ such that
$\bar{x}_i=e_{i}x_i$, for all $i$. The pair $(\aa,\bar{\;})$ is
called the {\it quantum torus with involution} determined by the
vector $\be$ and the matrix $\bq$. We have
$\aa=\aa_{+}\oplus\aa_{-}$ where
$\aa_{+}=\{h\in\aa\mid\bar{h}=h\}$ and
$\aa_{-}=\{s\in\aa\mid\bar{s}=-s\}$. Set
$$
I_\be=\{i\mid 1\leq i\leq\nu,\; e_i=-1\}\andd J_{\bq}=\{(i,j)\mid
1\leq i<j\leq\nu,\; q_{ij}=-1\}.
$$
Put
$$
\zeq=\{\d\in\zn\mid\sum_{i\in I_\be}n_i+\sum_{(i,j)\in J_{\bq}}
n_i n_j=0\}\andd\zeq^{c}=\zn\setminus\zeq.
$$
Then $\zeq$ is a semilattice in ${\mathbb R}^{\nu}$.
%Let $m\geq
%1$ and consider $m$ distinct coset representatives
%$\tau_1,\ldots,\tau_m$ of $2\zn$ in $\zn$.

Fix $n\geq 1$ and let $\mna$ be the Lie algebra of $n\times n$
matrices with entries from $\aa$.
%Make $\mna$ to a
%$\zn$--graded Lie algebra as follows. Set
%$$
%\lam_1=0,\ldots,\lam_{2\ell}=0,\quad\lam_{2\ell+1}=
%\tau_1,\ldots,\lam_{2\ell+m}=\lam_n=\tau_{m}.
%$$
%For $\d\in\zn$ and $1\leq i,j\leq n$, define
%$\deg(x^{\d}e_{ij})=2\d+\lam_{i}-\lam_{j}$.
Next define a bilinear form on $\mna$ as follows. Define $\ep$ in
the dual space of $\aa$ by the linear extension of $\ep(1)=1$ and
$\ep(\xd)=0$ for any nonzero $\d\in\zn$. Then for $A,B\in\mna$
define $(A,B)=\ep(\tr(AB))$. This defines a symmetric
nondegenerate invariant bilinear form on $\mna$.

Let $\kk$ be the Lie subalgebra $\slna$ of $\mna$ consisting of
matrices $A$ such that $\tr(A)\equiv 0$, $\mod [\aa,\aa]$.
%The
%$\zn$--grading on $\mna$ induces a $\zn$--grading on $\slna$ by
%$$\kk=\bigoplus_{\d\in\zn}\ks\quad\hbox{where}\quad\ks=\kk\cap\mna^{\d}.
%$$
Let $\hd$ be the abelian subalgebra of $\slna$ with basis
$e_{ii}-e_{i+1,i+1}$, for $1\leq i\leq n-1$. Define
$\ep_i\in\hd^{\star}$ by $\ep_i(e_{jj}-e_{j+1,j+1})=\d_{ij}$.
Then, with respect to $\hd$, we have the root space decomposition
$$\kk=\bigoplus_{\dot{\a}\in\rd}\kk_{\dot{\a}},\hbox{ where }\rd=
\{\pm(\ep_{i}-\ep_{j})\mid 1\leq i<j\leq n-1\}.
$$
Suppose that $\kk$ has a $\zn$ grading, say
$\kk=\sum_{\d\in\zn}\kd$, such that
\begin{equation}\label{setting2}
\begin{array}{ll}
\d,\tau\in\zn,\;\;\d+\tau\not=0\Longrightarrow(\kk^\d,\kk^{\tau})=
\{0\},\quad\hbox{and}
\vspace{2mm}\\
 \{\d\in\zn\mid\kd\not=\{0\}\}\hbox{ spans }{\bbbc}^\nu.
\end{array}
\end{equation}
Later we would like to consider the fixed point subalgebra of
$\kk$ under some automorphism $\sg$. According to our previous
notation, we should write $\ks$ for the corresponding fixed point
subalgebra. But this might cause some confusion with the notation
which we used for the grading on $\kk$. To prevent this, we
devote the upper index $\sg$ only to indicate the fixed points,
and we use other symbols such as $\d,\tau,\gamma,\ldots$ for the
grading on $\kk$.

For $1\leq i\leq\nu$ define $d_i\in Der(\kk)$ by $d_i(x)=n_ix$
for $x\in\kk^{\d}$, $\d=(n_1,\ldots,n_\nu)$. Set
$\dd=\sum_{i=1}^{\nu}{\bbbc}d_i$. Consider a $\nu$-dimensional
vector space $\cc=\sum_{i=1}^{\nu}{\bbbc}c_i$ and put
$$\gg=\kk\oplus\cc\oplus\dd\andd\hh=\hd\oplus\cc\oplus\dd.$$
Let $\{\d_1,\ldots,\d_\nu\}$ be the basis of $\dd^{\star}$ dual to
$\{d_1,\ldots,d_\nu\}$ and consider them as elements of
$\hh^{\star}$ by $\d_i(\hd\oplus\cc)=\{0\}$. Consider $\zn$ as a
subset of $\hh^{\star}$ through
$(n_1,\ldots,n_{\nu})=\sum_{i=1}^{\nu}n_{i}\d_{i}$.

Extend the bracket on $\kk$ to $\gg$ by,
$$[\cc,\cc]_{\gg}=\{0\}=[\dd,\dd]_\gg,\;\;[d_i,x]_\gg =
d_i(x),\;\;[x,y]_\gg=[x,y] +\sum_{i=1}^{\nu}(d_ix,y)c_i,$$ for
$x,y\in\kk$. Then $\gg$ is a Lie algebra. Also extend the form on
$\kk$ to $\gg$ by requiring $(\cc,\cc)=\{0\}$, $(\dd,\dd)=\{0\}$,
$(\cc\oplus\dd,\kk)=\{0\}$ and $(c_i,d_j)=\d_{ij}$. Then $\fm$ is
a symmetric, nondegenerate invariant bilinear form on $\gg$. It
is shown in [BGK] that $(\gg,\fm,\hh)$ is a tame EALA of type
$A_{n-1}$ and that $\gc=\kk\oplus\cc$.

In the next few examples we consider an automorphism $\sg$ of
$\kk$ ( or of a subalgebra of $\kk$) and we extended it to an
automorphism of $\gg$ by acting as identity on $\cc\oplus\dd$.
Then we have $\gs=\ks\oplus\cc\oplus\dd$ and
$\hs=\hds\oplus\cc\oplus\dd$, where $\gs$, $\ks$, $\hs$ and
$\hds$ denote the fixed points of $\gg$, $\kk$, $\hh$ and $\hd$
under $\sg$, respectively.

\begin{exa}\label{3-6}
Fix $\ell\geq 1$. In the previous setting let $\be=1_{\nu}$,
$\bq=1_{\nu\times\nu}$ and $n=2\ell$. Define $\deg(\xd
e_{pq})=2\d$ for $\xd\in\kd$. This defines a $\zn$ grading on
$\kk$ satisfying (\ref{setting2}). Define $\sg\in\hbox{Aut}(\kk)$
by
$$\sg(Y)=-KY^{t}K\quad\hbox{where}\quad
K=\begin{pmatrix} 0&I_{\ell}\\
I_{\ell}&0
\end{pmatrix}.
$$
$\sg$ is a period 2 automorphism of $(\gg,\fm,\hh)$, so it
satisfies A1. It is easy to see that $\sg$ satisfies conditions A2
and A3. To see A4 holds, note that the root system of $\gg$ with
respect to $\hh$ is
$$R=2\zn\cup\{\ep_i-\ep_j+2\d\mid\d\in\zn,\;\;1\leq i\not=j\leq
n-1\},
$$ where $\ep_i$ is the element of $\hh^{\star}$ which acts on
$\hh$ by $\ep_i(e_{jj}-e_{j+1,j+1})=\d_{ij}$. Also note that
$\hds=\{\sum_{i=1}^{\ell}a_i(e_{ii}-e_{\ell+i,\ell+i})\mid
a_i\in{\bbbc}\}$. Then it is clear that $$0\not=\a\in R\Rightarrow
\pia=\a_{|_{\hs}}\not=0.$$ Thus, by [ABP, Proposition 3.25], A4
holds. One can see that
$$
\gs=\{\begin{pmatrix}A&B\\C&-A^{t}\end{pmatrix}\mid
B^{t}=-B,\;C^{t}=-C,\;A,B,C\in M_{\ell}(\aa)\}.
$$
The root system of $\gs$ with respect to $\hs$ is
$$\rs=2\zn\cup\{\pm(\ep_{i}\pm\ep_{j})+2\d\mid\d\in\zn,\;\;1\leq
i\not=j\leq\ell\},
$$
which is an irreducible reduced EARS of type $D_{\ell}$. Thus by
Corollary \ref{2-28}, $\gs$ is an EALA of type $D_{\ell}$. It is
not difficult to show that $\gsc=\ks\oplus\cc$. Thus
$\gs\cap\gc=\gsc$. Therefore by Lemma \ref{2-16b} $\gs$ is a tame
EALA.
\end{exa}

\begin{exa}\label{3-7}
Let $(\aa,\bar{\;})$ be the quantum torus determined by $\be$ and
$\bq$ as in (\ref{setting1}). Suppose that $m\geq 1$ and let
$\tau_1,\ldots\tau_m$ represent distinct cosets of $2\zn$ in
$\zn$ with $\tau_1=0$ and $\tau_i\in\zeq$ for all $i$. Define
$\deg(\xd e_{pq})=2\d+\lam_{p}-\lam_{q}$ where $\lam_1=\cdots =
\lam_{2\ell}=0$ and $\lam_{2\ell+i}=\tau_{i}$ for $1\leq i\leq m$.
Set
$$ F=\begin{pmatrix}
x^{\tau_1}&\hdots&0\\
\vdots&\ddots&\vdots\\
0&\hdots&x^{\tau_m}\\
\end{pmatrix}
\andd G=\begin{pmatrix}
0&I_{\ell}&0\\
I_{\ell}&0&0\\
0&0&F\end{pmatrix},
$$
where $I_{\ell}$ is the $\ell\times\ell$ identity matrix. Then
$G$ is a $n\times n$ matrix with $n=2\ell+m$. For $Y\in\kk=\slna$
define $\sg(Y)=-G^{-1}\bar{Y}^{t}G$. Then $\sg$ defines a period
two automorphism of $\gg$. It is straightforward to see that
$\sg$ satisfies A2--A3. Elements of $\ks$ are of the form
$$
X=\begin{pmatrix} A&S&-\bar{D}^{t}F\\
T&-\bar{A}^{t}&-\bar{C}^{t}F\\
C&D&B
\end{pmatrix},
$$
where $\bar{S}^{t}=-S$, $\bar{T}^{t}=-T$,
$F^{-1}\bar{B}^{t}F=-B$, and $\hbox{tr}(X)\equiv 0$ $\hbox{mod
}[\aa,\aa]$.

We next want to check A4. Note that
$\hds=\{\sum_{i=1}^{\ell}a_i(e_{ii}-e_{\ell+i,\ell+i})\mid
a_i\in\mathbb C\}$. It is easy to see that
$$\gg=\bigoplus_{\a\in\hh^{\star}}\gg_{\a}=
\sum_{\d\in\zn}\sum_{\a\in\rd}\gg_{\a+\d},$$ where
$\rd=\{\pm(\ep_i-\ep_j)\mid 1\leq i\leq j\leq n\}$,
$$\gg_{0}=\hd\oplus\cc\oplus\dd,\andd
\gg_{\a+\d}=\kk_{\a}\cap\kk^{\d}\;\hbox{ if }\;\a+\d\not=0.$$ In
particular, the root system $R$ of $\gg$ is of the form
$$R=\{\a+\d\mid\a\in\rd,\;\d\in\zn\hbox{ and
}\kk_{\a}\cap\kk^{\d}\not=\{0\}\}.
$$
To describe the root system of $\gg$ note that if $1\leq i<j\leq
2\ell$, then $\deg(e_{ij})=\lam_i-\lam_j=0$ and so
$e_{ij}\in\kk_{\ep_i-\ep_j}\cap\kk^{0}$. Therefore
$\ep_i-\ep_j\in R$. If $1\leq i\leq 2\ell$ and $1\leq j\leq m$,
then $\deg(\xd e_{i,2\ell+j})= 2\d+\lam_i-\lam_{j}=2\d-\tau_j$.
Therefore $\xd e_{i,2\ell+j}\in\kk_{\ep_{i}-\ep_{2\ell+j}}
\cap\kk^{2\d-\tau_j}$. So $\ep_i-\ep_{2\ell+j}+2\d-\tau_j\in R$.
Finally if $1\leq i\leq j\leq m$, then $\deg(\xd
e_{2\ell+i,2\ell+j})=2\d+\tau_i-\tau_j$. So
$$\xd e_{2\ell+i,2\ell+j}\in\kk_{\ep_{2\ell+i}-
\ep_{2\ell+j}}\cap\kk^{2\d+\tau_i-\tau_j}.$$ Thus
$\ep_{2\ell+i}-\ep_{2\ell+j}+2\d+\tau_i-\tau_j\in R$. This in
particular shows that $\ep_{2\ell+i}-\ep_{2\ell+j}\not\in R$, as
$\tau_i-\tau_j\not\in 2\zn$ if $i\not= j$. Set
$$
\ep'_{i}=\ep_{i}\hbox{ for }1\leq i\leq 2\ell\andd
\ep'_{2\ell+i}=\ep_{2\ell+i}+\tau_i\hbox{ for }1\leq i\leq m.
$$
Then
$$
R=2\zn\cup\{\pm(\ep'_{i}-\ep'_{j})+2\d\mid \d\in\zn,\;\; 1\leq
i\leq j\leq n\}.
$$
It now follows that
$$
\a\in R\setminus\{0\}\Rightarrow\pia=\a_{|_{\hs}}\not=0.
$$
Thus by [ABP, Proposition 3.25], the condition A4 also holds. It
is shown in [A2, Proposition 3.23], the root system $\rs$ of $\gs$
with respect to $\hs$ is an irreducible reduced EARS of type
$$\begin{array}{ll}
A_1&\hbox{if }\ell=1,\;\be=1_{\nu}\hbox{ and }\bq=1_{\nu\times\nu}\\
B_{\ell}&\hbox{if }\ell\geq 2,\;\be=1_{\nu}\hbox{ and }\bq=1_{\nu\times\nu}\\
BC_{\ell}&\hbox{if }\be\not=1_{\nu}\hbox{ or }\bq\not=1_{\nu\times\nu}.\\
\end{array}
$$
Thus by Corollary \ref{2-28}, $\gs$ is an EALA of the above type.
It is not difficult to show that $\gsc=\ks\oplus\cc$. Since
$\gc=\kk\oplus\cc$, we have $\gsc=\gs\cap\gc$. It now follows
from Lemma \ref{2-16b} that $\gs$ is a tame EALA.
\end{exa}

\begin{exa}\label{3-8}
Suppose $\ell\geq 2$ and $n=2\ell$. Let $\kk=\slna$ and consider
the EALA $(\gg,\fm,\hh)$ as in the setting before Example
\ref{3-6}. Define a gradation on $\kk$ by $\deg(\xd e_{pq})=2\d$.
This gradation satisfies (\ref{setting2}). Let $K=\begin{pmatrix}
0&I_{\ell}\\
-I_{\ell}&0
\end{pmatrix}
$ and for $Y\in\kk$ define $\sg(Y)=-K^{-1}\bar{Y}^{t}K$. Extend
$\sg$ to a period 2 automorphism of $\gg$. Check that A2--A3 hold
for $\sg$. Let $X\in\kk$. It is easy to see that $X\in\ks$ if and
only if
\begin{equation}\label{3-8a}
X=\begin{pmatrix}
A&S\\T&-\bar{A}^{t}\end{pmatrix}\quad\hbox{with}\quad\bar{S}^{t}=
S\andd\bar{T}^{t}=T.
\end{equation}
We also have
$$\hds=\{\sum_{i=1}^{\ell}a_{i}(e_{ii}-e_{\ell+i,\ell+i})\mid
a_{i}\in\bbbc\}.$$ Next we check A4. Note that the root system
$R$ of $\gg$ is of the form
$$
R=2\zn\cup\{\ep_{i}-\ep_{j}+2\zn\mid 1\leq i\not=j\leq n\},
$$
(see Example \ref{3-6}). It is now easy to see that
$$
\a\in R\setminus\{0\}\Rightarrow \pia=\a_{|_{\hs}}\not=0.
$$
Thus by [ABP,Proposition 3.25], $\sg$ satisfies A4. The root
system $\rs$ of $\gs$ with respect to $\hs$ is an irreducible
reduced EARS of type $C_{\ell}$ (see [AABGP, III.Theorem 4.7]).
Thus by Corollary \ref{2-28}, $\gs$ is an EARS of type $C_{\ell}$.
As in the previous examples we can see that $\gs$ is a tame EALA.
\end{exa}

\begin{exa}\label{3-9}
Let $\nu\geq 1$ and let $\aa$ be the associative commutative
Laurent polynomials in $\nu$--variables $x_1,\ldots,x_\nu$. Let
$\ell\geq 1$, $m\geq 1$, $n=\ell+m$. Let
$$
\kk=\{\begin{pmatrix}A&B\\
C&-A^{t}\end{pmatrix}\mid B^t=-B,\;C^t=-C,\;A,B,C\in \mna\}$$ and
$$\hd=\{\sum_{i=1}^{n}a_i(e_{ii}-e_{n+i,n+i})\mid
a_i\in\bbbc\}.$$ Let $\gg=\kk\oplus\cc\oplus\dd$ and
$\hh=\hd\oplus\cc\oplus\dd$. Then by Example \ref{3-6},
$(\gg,\fm,\hh)$ is an EALA of type $D_{n}$ $(n\geq 4)$.

Let $\tau_1,\ldots\tau_m$ represent distinct cosets of $2\zn$ in
$\zn$ with $\tau_1=0$. Define $\deg(\xd
e_{pq})=2\d+\lam_{p}-\lam_{q}$ where
$$\lam_i=\lam_{n+i}=0\hbox{ for }1\leq i\leq\ell\andd
\lam_{\ell+i}=-\tau_i,\;\;\lam_{n+\ell+i}=\tau_i\hbox{ for }
1\leq i\leq m.
$$
This defines a gradation on $\kk$ satisfying (\ref{setting2}).

Put
$$
F=\begin{pmatrix}x^{\tau_{1}}&\cdots&0\\
\vdots&\ddots&\vdots\\
0&\cdots&x^{\tau_{m}}\end{pmatrix}\andd K=\begin{pmatrix}
I_{\ell}&0&0&0\\
0&0&0&F\\
0&0&I_{\ell}&0\\
0&F^{-1}&0&0
\end{pmatrix}.
$$
Define a period 2 automorphism of $\kk$ by $\sg(Y)=KYK$ and extend
it to an automorphism of $\gg$ as before. It is easy to check that
$\sg$ satisfies A2--A3. Next we check A4. It is not difficult to
see that the root system $R$ of $\gg$ is of the form
$$R=2\zn\cup\{\pm(\ep'_i\pm\ep'_j)+2\zn\mid 1\leq i<j\leq n\},$$
where
$$
\ep'_i=\ep_i\hbox{ for }1\leq i\leq
\ell\andd\ep'_i=\ep_{\ell+i}-\tau_{i}\hbox{ for }1\leq i\leq m.
$$
Also we have
$$
\hs=\{\sum_{i=1}^{\ell}a_i(e_{ii}-e_{n+i,n+i})\mid
a_i\in\bbbc\}\oplus\cc\oplus\dd.
$$
It follows that
$$
\a\in R\setminus\{0\}\Rightarrow \pia=\a_{|_{\hs}}\not=0.
$$
This proves that A4 holds for $\sg$. One can check that each
element of $\ks$ is of the form
\begin{equation}\label{3-10}
\begin{pmatrix}A&-D^t&S&-D^tF\\
FC&-B^t&FD&FPF\\
T&-C^t&-A^t&-C^tF\\
C&P&D&B
\end{pmatrix}\quad\hbox{with}\begin{array}{c}
S^t=-S,\;T^t=-T,\;P^t=-P\\
F^{-1}B^tF=-B,
\end{array}
\end{equation}
where $A,S,T\in M_{\ell}(\aa)$, $C,D\in M_{m\times\ell}(\aa)$ and
$B,P\in M_{m}(\aa)$. It is shown in [A2, Proposition 3.24] that
the root system $\rs$ of $\gs$ with respect to $\hs$ is an
irreducible reduced EARS of type
$$
A_{1}\quad\hbox{if }\ell=1\andd B_{\ell}\quad\hbox{if }\ell\geq 2.
$$
Thus by Corollary \ref{2-28}, $\gs$ is an EALA of the above type.

Note that $\gg$ is a tame EALA and $\gc=\kk\oplus\cc$ (see
[AABGP, Chapter III]). Let $\ks_{c}$ be the subalgebra of $\ks$
consisting of all matrices of the form (\ref{3-10}) with
$P=BF^{-1}$. By [A2, Lemma 3.17], $\gsc=\ks_{c}\oplus\cc$. It
follows that $\gsc=\gs\cap\gc$ if and only if $m=1$. Now Lemma
\ref{2-16b} implies that $\gs$ is tame if $m=1$. In fact one can
see that $\gs$ is tame if and only if $m=1$ (see [A2, Corollary
3.22]). This example in particular shows that $\gsc$ in general is
a proper subalgebra of $\gs\cap\gc$.
\end{exa}

\begin{exa}\label{3-11}
Consider the quantum torus $(\aa,\bar{\;})$ determined by the
vector $\be$ and the matrix $\bq$. Assume that $\be\not=1_{\nu}$
or $\bq\not=1_{\nu\times\nu}$. Then $\zeq^{c}\not=\emptyset$. Let
$\ell\geq 1$ and $m\geq 1$. Set $n:=\ell+m$. Consider the
subalgebra $\kk$ of $sl_{2n}(\aa)$ consisting of matrices of the
form (\ref{3-8a}). Let
$\hd=\{\sum_{i=1}^{n}a_i(e_{ii}-e_{n+i,n+i})\mid a_i\in\bbbc\}$.
Set $\gg=\kk\oplus\cc\oplus\dd$ and $\hh=\hd\oplus\cc\oplus\dd$.
By Example \ref{3-8}, $(\gg,\fm,\hh)$ is a tame EALA of type
$C_{\ell}$. Let $\tau_1,\ldots,\tau_m$ represent distinct cosets
of $2\zn$ in $\zn$ such that $\tau_i\in\zeq^{c}$. Define a
gradation on $\kk$ by $\deg(e_{pq})=2\d+\lam_{p}-\lam_{q}$ where
$\lam_{i}$'s are defined as in Example \ref{3-9}. Put
$$ F=\begin{pmatrix}
x^{\tau_1}&\hdots&0\\
\vdots&\ddots&\vdots\\
0&\hdots&x^{\tau_m}\\
\end{pmatrix}
\andd K=\begin{pmatrix}
0&0&I_{\ell}&0\\
0&-F^{-1}&0&0\\
-I_{\ell}&0&0&0\\
0&0&0&F\end{pmatrix}.
$$
Define a period 2 automorphism of $\kk$ by
$\sg(Y)=-K^{-1}\bar{Y}^{t}K$ and extend it to an automorphism of
$\gg$ as before. We have seen in the previous examples that $\sg$
satisfies A3. It is easy to check that $\sg$ satisfies A2. The
root system $R$ of $\gg$ is of the form
$$
R=2\zn\cup\{\pm(\ep'_i\pm\ep'_j)\mid 1\leq i\not= j\leq
n\}\cup\{2\ep'_i\mid 1\leq i\leq n\},
$$
where
$$ \ep'_i=\ep_i\;\;\hbox{ for }1\leq i\leq\ell\andd
\ep'_{\ell+i}=\ep_{\ell+i}-\tau_{i}\;\;\hbox{ for }1\leq i\leq m.
$$
Also note that
$\hds=\{\sum_{i=1}^{\ell}a_{i}(e_{ii}-e_{n+i,n+i})\mid
a_i\in\bbbc\}$. It follows that
$$
\a\in R\setminus\{0\}\Rightarrow\pia=\a_{|_{\hs}}\not=0.
$$
Thus $\sg$ satisfies A4. It is shown in [A2, (3.27)] that
$X\in\ks$ if and only if
\begin{equation}\label{3-11a}
X=\begin{pmatrix} A&-\bar{D}^t&S&\bar{D}^t\\
-FC&-\bar{B}^t&-FD&FPF\\
T&\bar{C}^t&-\bar{A}^t&-\bar{C}^tF\\
C&P&D&B
\end{pmatrix}\;\;\hbox{ with }\;\;\begin{array}{l}
\bar{S}^t=S,\;\;\bar{T}^t=T,\;\;\bar{P}^t=P,\\
F^{-1}\bar{B}^tF=-B,\hbox{ and}\\
\tr(X)\equiv 0\;\;\mod [\aa,\aa],\end{array}
\end{equation}
where $A,S,T\in M_{\ell}(\aa)$, $C,D\in M_{m\times\ell}(\aa)$ and
$B\in M_{m}(\aa)$. It is also shown in [A2, Proposition 3.50]
that the root system of $\gs$ with respect to $\hs$ is an
irreducible reduced EARS of type $BC_{\ell}$ $(\ell\geq 1)$. Thus
$(\gs,\fm,\hs)$ is an EALA of type $BC_{\ell}$

The Lie algebra $\gs$ is not tame. To see this let $\ks_c$ be the
subalgebra of $\ks$ consisting of matrices of the form
(\ref{3-11a}) with $P=-BF^{-1}$. By [A2, Lemma 3.54] if the
intersection of $\aa_{-}$ and the center of $\aa$ is nonzero,
then $\gs$ is not tame. Also if this intersection is zero, then
the matrix $X$ of the form (\ref{3-11a}) with $A=S=T=0$, $C=D=0$,
$P=I_m$ and $B=F$ is orthogonal to $\gsc$ but $X\not\in\gsc$.
Thus $\gs$ is not tame. It follows from Lemma \ref{2-16b} that
$\gsc$ is a proper subspace of $\gs\cap\gc$
\end{exa}

\begin{exa}\label{3-12} (See [ABP])
Let $(\gg,\fm,\hh)$ be an EALA with root system $R$. Consider the
so called affinization
$$\aff(\gg)=(\gg\otimes{\bbbc}[t,t^{-1}])\oplus{\bbbc}c\oplus{\bbbc}d,
$$
of $\gg$ introduced in [ABP]. Then $\aff(\gg)$ is a Lie algebra
where $c$ is central, $d=t\frac{d}{dt}$ is the degree derivation
so that $[d,x\otimes t^n]=nx\otimes t^n$, and
$$ [x\otimes t^n,y\otimes t^m]=[x,y]\otimes
t^{n+m}+n(x,y)\d_{m+n,0}c. $$ Extend the form $\fm$ to
$\aff(\gg)$ such that $c$ and $d$ are dually paired. Then the
triple
$$\big(\aff(\gg),\fm,\hh\oplus{\bbbc}c\oplus{\bbbc}d\big)$$
is an EALA with root system $\rt=R+{\bbbz}\d$ where
$\d\in\hht^{\star}$ is defined by $\d(d)=1$ and
$\d(\hh+{\bbbc}c)=0$. Moreover, $\aff(\gg)$ is tame if and only if
$\gg$ is tame. Next consider an automorphism $\sg$ of $\gg$
satisfying A1-A4. Extend $\sg$ to an automorphism of $\aff(\gg)$
by
$$\sg(x\otimes t^i+rc+sd)=\zeta^{-i}\sg(x)\otimes t^i+rc+sd.
$$
It is easy to see that $\sg$ as an automorphism of $\aff(\gg)$
satisfies A1-A4. By [ABP, Lemmas 3.49, and 3.62] the root system
of $\aff(\gg)^{\sg}$ is irreducible. Thus by Corollary \ref{2-28},
$\aff(\gg)^{\sg}$ is an EALA. Therefore all Lie algebras $\gt$
constructed in [ABP] can be considered as fixed point subalgebras
of the loop algebra of some EALA.
\end{exa}

 \scriptsize{\tiny DEPARTMENT OF MATHEMATICS, UNIVERSITY
OF ISFAHAN, ISFAHAN, IRAN, P.O.BOX 81745.

DEPARTMENT OF MATHEMATICS AND STATISTICS,  UNIVERSITY OF
SASKATCHEWAN, SASKATOON, SASKATCHEWAN, CANADA, S7N 5E6.

DEPARTMENT OF MATHEMATICS, UNIVERSITY OF ISFAHAN, ISFAHAN, IRAN,
P.O.BOX 81745.}

\end{document}